%% file: fastigabem.tex
\pdfoutput=1
\documentclass[3p]{elspreprint}

\input{settings}

\input{commands}

\begin{document}

\input{\CommonPath/title}


\input{\CommonPath/introduction}

\input{\CommonPath/bie}

\input{\CommonPath/igabem}

\input{\CommonPath/hmat}

\input{\CommonPath/results}

\input{\CommonPath/conclusion}

\clearpage

\appendix

\input{\CommonPath/appendixDerivatives}

\input{\CommonPath/appendixOperations}

\input{\CommonPath/integration}

\input{\CommonPath/literature}

\end{document}

%% file: settings.tex
\usepackage[intlimits]{amsmath}
\usepackage{amssymb} 
\usepackage{amsthm}
\usepackage{amsfonts}
\usepackage{bm}
\usepackage{mathdots} 
\usepackage{siunitx} 

\usepackage[utf8]{inputenc} 
\usepackage[english]{babel}  
\usepackage[T1]{fontenc}

\usepackage{times}
\usepackage{mathptmx}  
\DeclareSymbolFont{cmletters}{OML}{cmr}{m}{n}
\DeclareMathSymbol{\epsilon}{0}{cmletters}{'017} 
\DeclareMathAlphabet{\mathcal}{OMS}{cmsy}{m}{n} 

\usepackage{booktabs} 
\usepackage{hyperref} 


%% file: commands.tex
\newcommand{\rk}{\mathcal{R}_k} 
\newcommand{\h}{\mathcal{H}} 
\newcommand{\hmatrix}{\mbox{$\h$-matrix}}
\newcommand{\hmatrices}{\mbox{$\h$-matrices}}

\newcommand{\orderof}[1]{\mathcal{O}\left(#1\right)} 
\newcommand{\order}{\mathcal{O}} 
\newcommand{\storage}{\mathsf{St}} 
\newcommand{\operations}{\mathsf{Op}} 

\newcommand{\R}{\mathbb{R}}
\newcommand{\N}{\mathbb{N}}

\newcommand\op[1]{\mathcal{#1}}
\newcommand\operate[1]{(#1)}

\providecommand\url[1]{\emph{#1}}

\DeclareMathOperator{\dx}{d}
\DeclareMathOperator{\Tr}{Tr}

\DeclareMathOperator{\supp}{supp}

\newcommand\fromto[2]{\{ #1, \dots, #2 \}}

\newcommand\uu{r}
\newcommand\vv{s}

\newcommand\mIndexI{\vek{i}}
\newcommand\mIndexJ{\vek{j}}
\newcommand\mIndexP{\vek{p}}

\newcommand\greville{\bar{\uu}}

\renewcommand{\div}{\nabla \cdot}
\newcommand{\grad}{\nabla}
\newcommand{\curl}{\nabla \times}

\newcommand\dgamma[1]{\dx \mathrm{s}_{\pt{#1}}}

\newcommand{\dTr}{\op{T}} 

\newcommand\vek[1]{\mathbf{#1}}
\newcommand\veks[1]{\boldsymbol{#1}}
\newcommand\mat[1]{\mathbf{#1}}
\newcommand\tens[2]{\mathsf{#1}}

\newcommand{\trans}{\intercal} 
\newcommand{\inv}{{-1}} 

\newcommand\utens{u} 
\newcommand\ttens{t} 

\newcommand\uvek{\vek{u}} 
\newcommand\tvek{\vek{t}} 
\newcommand\uvekc{\tilde{\uvek}} 
\newcommand\tvekc{\tilde{\tvek}} 

\newcommand\ct{T} 
\newcommand\bct{T_{I\times J}} 

\newcommand\fund[1]{\tens{#1}{2}}

\newcommand\be{\tau} 
\newcommand\ieg{\tilde{\tau}} 
\newcommand\ieKs{\hat{\tau}} 
\newcommand\ieAr{\hat{\tau}^{\boxplus}} 
\newcommand\ie{\hat{\tau}^{\boxdot}} 

\newcommand\beistrich{,\,}
\newcommand\und{\,\text{and}\,}

\newcommand\with{\,\text{with}\,}

\newcommand\pt[1]{\veks{#1}}
\newcommand\ofpt[1]{(\pt{#1})}

\newcommand\eps{\varepsilon}
\newcommand\err{\epsilon}
\newcommand\mychi{\mathcal{X}}

\newcommand{\tikzfig}[5]{
  \begin{figure}[ht]
    \centering
    \includegraphics{#5}
    \caption{#3}
    \label{fig:#4}
  \end{figure}
}

\newenvironment{mytable}[3]%
{
  \begin{table}[ht]
    \def\mycap{#1}
    \def\mylabel{tab:#2}
    \centering
    \begin{tabular}{#3}
      \toprule
}
{
  \bottomrule
  \end{tabular}
  \caption{\mycap}
  \label{\mylabel}
  \end{table}
}
\newcommand{\mytableheader}[1]{
       #1 \\ \midrule
}

\newcommand{\mycomment}[2][-]{ }%
\newcommand{\todoref}[2][reference]{}

\newcommand{\figref}[1]{Figure~\ref{#1}}
\newcommand{\Figref}[1]{Figure~\ref{#1}}
\newcommand{\secref}[1]{section~\ref{#1}}

\newcommand{\Tabref}[1]{Table~\ref{#1}}

\newcommand{\CommonPath}{.}

\newcommand\revdel[2]{}
\newcommand\revcomment[1]{} 
\newcommand\revadd[2]{#1}
\newcommand\revmod[2]{#1}


%% file: title.tex
\title{Fast Isogeometric Boundary Element Method based on Independent Field Approximation}

\begin{frontmatter}

\author[ifbaddr]{Benjamin Marussig}

\author[ifbaddr]{Jürgen Zechner\corref{cor1}}

\author[ifbaddr,newcastleaddr]{Gernot Beer}

\author[ifbaddr]{Thomas-Peter Fries}

\address[ifbaddr]{Institute of Structural Analysis, Graz University
  of Technology, Lessingstraße 25/II, 8010 Graz, Austria}

\address[newcastleaddr]{Centre for Geotechnical and Materials Modelling, University of Newcastle,
  Callaghan, NSW 2308, Australia}

\cortext[cor1]{Corresponding author. 
  Tel.: +43 316 873 6181, fax: +43 316 873 6185, mail: \url{ifb@tugraz.at}, web: \url{www.ifb.tugraz.at}}

\input{\CommonPath/abstract}

\end{frontmatter}


%% file: abstract.tex
\begin{abstract}
An isogeometric boundary element method for problems in elasticity is presented, which is based on an independent approximation for the geometry, traction and displacement field.
 This enables a flexible choice of refinement strategies, permits an efficient evaluation of geometry related information, a mixed collocation scheme which deals with discontinuous tractions along non-smooth boundaries and a significant reduction of the right hand side of the system of equations for common boundary conditions.
 All these benefits are achieved without any loss of accuracy compared to conventional isogeometric formulations.
\revdel{A tailored integration procedure maintains the precision for the computation of the coefficients of the system matrices. These matrices}{Removed because the integration section has been removed from the main text} 
\revadd{The system matrices}{} are approximated by means of hierarchical matrices to reduce the computational complexity for large scale analysis. For the required geometrical bisection of the domain, a strategy for the evaluation of bounding boxes containing the supports of NURBS basis functions is presented. The versatility and accuracy of the proposed methodology is demonstrated by convergence studies showing optimal rates and real world examples in two and three dimensions.
\end{abstract}



%% file: introduction.tex
\section{Introduction}
\label{sec:introduction}

Isogeometric analysis (IGA) \cite{Hughes2005a} has received much attention in recent years and it has been successfully applied to a variety of applications.
It offers precise and efficient geometric modeling, refinement without the need for communication with the design model and control over the smoothness of the basis functions.  
The most important feature of IGA is that the discretisation is based on the computer aided design (CAD) representation and therefore it is geometrically exact regardless of the fineness of the underlying basis.
This unified representation facilitates the integration of design and analysis model.
Hence the overall analysis time, including the creation and improvement of the analysis model, can be reduced significantly \cite{Cottrell2009b}.

However, the main challenge in the context of an isogeometric finite element analysis (FEA) is the generation of a volume discretisation from CAD models which are usually based on a boundary representation (B-Rep). 
This task is far from trivial and is still an open research topic \cite{Thurston1982a,Thurston1997b}. 
For linear problems, a volume discretisation can be completely avoided by using an analysis which is based on boundary integral equations, such as the boundary element method (BEM).
The distinguishing feature of BEM is that it is based on a B-Rep as well. Hence the CAD model can be used for the simulation without the need for a volume discretisation.

This natural connection is possibly the reason for some early attempts to include CAD representations in a BEM framework. \revadd{For example, spline collocation methods were used to develop convergence estimates for the BEM in two~\protect\cite{Arnold1983a,Arnold1985a,Saranen1985a,Costabel1987a,Saranen1988a}
and three dimensions~\protect\cite{Arnold1984a,Schneider1991a,Costabel1992a}.}{References have been added.}
A BEM formulation based on cubic splines was
\revdel{, for example,}{}%
proposed by \cite{Ligget1981a} and \cite{Yu1994a} to solve groundwater problems and the Laplace's equation, respectively.
Cubic B-splines were applied to potential problems \cite{Cabral1990a,Cabral1991a} and to coplanar waveguides \cite{Qi1993p}.
The analysis of electromagnetic problems with three dimensional B\'{e}zier patches has been demonstrated in \cite{Schlemmer1994p,Ushatov1994a}. \citet{Turco1998a} suggested to use three dimensional B-spline elements to solve elasticity problems.
In \cite{Maniar1995phd} B-spline surfaces are used by means of the panel method which is a variation of the BEM for solving the potential flow around aerodynamic bodies.
To the best of the authors' knowledge, the integration of rational non-uniform B-splines (NURBS) surfaces into a BEM analysis has been first investigated by the authors of \cite{Valle1994a} and \cite{Rivas1996a} in the context of the method of moments, a Galerkin approach to solve Maxwell's equations.

The introduction of IGA for FEA \cite{Hughes2005a,Cottrell2009b} has reignited the usage of CAD representation in the BEM community. 
The isogeometric BEM has been applied to the Laplace's equation \cite{Politis2009p,Gu2012a}, linear elasticity \cite{Li2011a,Simpson2012a,Beer2013a,Scott2013a,Marussig2014a}, Stokes flow \cite{Heltai2014a},  low-frequency acoustic problems \cite{Simpson2014a}, electro-magnetic problems \cite{Vazquez2012a} and an extended formulation for the Helmholtz equation in two dimensions \cite{Peake2013a}.

Despite the fact that IGA and BEM complement each other extremely well, the BEM naturally comes with a high degree of complexity: the system matrix is fully populated and the solution of the linear system requires a significant numerical effort. For the past three decades, much research has been dedicated to the reduction of the algorithmic complexity. The most important approaches are the fast multipole method (FMM) \cite{rokhlin1985}, hierarchical matrices (\hmatrices{}) \cite{hackbusch1999}, the wavelet method \cite{beylkin1991} and fast Fourier transformation based methods \cite{phillips1994}. 

Only a few publications deal with the application of fast BEM formulations to analysis involving B\'{e}zier, \revmod{B-spline}{removed spelling mistake} or NURBS functions. \citet{Gonzalez2006p} presented the combination of a multilevel FMM and B\'{e}zier patches. This is the first application of a fast solution method in that context. \citet{Harbrecht2010a} applied the wavelet method to trimmed NURBS patches by means of a decomposition into four-sided patches. In \cite{Harbrecht2013a} the authors compared different customised fast solution techniques on NURBS patches solving Laplace's problem. The same problem applied to the FMM in two dimensions has been studied in \cite{Takahashi2012a}.

\revdel{In this paper the first application of a fast solution technique to isogeometric BEM in elasticity is presented. The \hmatrices{} approach is employed with a scheme for the geometric bisection which is tailored towards NURBS patches.
}{Removed to avoid the impression that it is a paper about fast BEM only}%
\revadd{In this paper an efficient isogeometric BEM formulation for elasticity is introduced which is based on an independent discretisation for geometry, displacement field and traction field as well as the application of \hmatrices{}: }{New introduction describing the features of the paper}%
\revdel{In addition, t}{This sentence has been moved.}%
The proposed formulation allows discontinuous tractions along non-smooth boundaries while maintaining the physical constraint of continuous displacements.
\revdel{Due to this flexibility, the r}{This sentence has been moved.}%
Refinement can be specifically performed for fields where a finer resolution is necessary, while for others a coarse representation is maintained.
This results in an efficient evaluation of geometry related information and a minimisation of the right hand side of the system of equations and its storage requirements.
\revadd{The \hmatrix{} format is applied to the system matrices in order to reduce the computational complexity, in particular for the left hand side of the system of equations. In that context, the distinct modifications for the \hmatrix{} construction with NURBS patches are presented.}{This shall clarify that the paper does not present a new H-matrices technique.}
The proposed approach violates the \emph{isoparametric} concept. However, it is still \emph{isogeometric} since the same functions as in the CAD model are used and, hence the geometry representation is still exact.


The paper is organised in the following sections. First a brief overview on the boundary value problem and its corresponding boundary integral equation is given. The main part is found in \secref{sec:IGABEM}. Here the customised and individual discretisation with NURBS basis functions is described as well as the collocation scheme and the block-system of equations.
\revdel{Section 4 is devoted to the numerical integration of the resulting boundary integrals with an adaptive scheme, customised for NURBS patches. This is a crucial point to preserve accuracy.}{Removed because the integration section has been removed from the main text}%
In \secref{sec:hmat} the concept of \hmatrices{} and their application to the formulation is described, focusing on the geometric bisection of NURBS patches. Finally, numerical results demonstrate convergence and efficiency of the method presented.


%% file: bie.tex
\section{Boundary Integral Equation}
\label{sec:BIE}

An elastic body $\Omega$ subject to external loading without body forces is considered. Its behaviour in terms of displacements $\utens$ is described by the partial differential equation
\begin{align}
\label{eq:PDE}
  \op{L}\utens\ofpt{x} = -\left( \lambda + 2 \mu \right) 
  \div \grad \utens \ofpt{x} + \mu \curl (\curl \utens \ofpt{x}) =
  0 & & \pt{x} \in \Omega
\end{align}
where $\op{L}$ denotes the Lam\'{e}-Navier operator with the Lam\'{e} constants $\lambda$ and $\mu$ \cite{kupradze1979}. The closed boundary of the domain is denoted by $\Gamma$ and the surface normal $\vek{n}$ points outside the considered domain.
In the following sections,
\revdel{almost everything}{}%
\revadd{equation~\protect\eqref{eq:PDE} is modified, so that it}{}
is described by means of boundary data. Therefore the boundary trace
\begin{align}
  \label{eq:trace}
  \Tr \utens \ofpt{x} &= \lim_{\vek{x}\rightarrow\vek{y}} \utens \ofpt{x} =
  \utens \ofpt{y} & & \pt{x} \in \Omega,\, \pt{y} \in \Gamma
\end{align}
and the conormal derivative
\begin{align}
  \label{eq:conormal-derivative}
  \dTr \utens \ofpt{x} & = \lambda \div \utens \ofpt{y} \vek{n}\ofpt{y} +
  2\mu \grad \utens \ofpt{y} \cdot \vek{n}\ofpt{y} + \mu \vek{n}\ofpt{y}
  \times (\curl \utens \ofpt{y}) & & \pt{x} \in \Omega,\, \pt{y} \in
  \Gamma
\end{align}
are introduced \cite{steinbach2008}.
The trace operator \eqref{eq:trace} maps displacements~$\utens\ofpt{x}$ in the domain to boundary displacements~$\utens\ofpt{y}$. With the material law, the conormal derivative \eqref{eq:conormal-derivative} maps $\utens\ofpt{x}$ to surface traction $\ttens\ofpt{y}$. 

The boundary is split into a 
\revmod{Neumann part $\Gamma_N$ and a Dirichlet part $\Gamma_D$, such that }{}%
${\Gamma = \Gamma_N \cup \Gamma_D}$ and ${\Gamma_N \cap \Gamma_D = \emptyset}$. This leads to the following homogeneous boundary value problem (BVP): %
Find a displacement field $\utens\ofpt{x}$ so that
\begin{equation}
\label{eq:BVP}
  \begin{aligned}
    \op{L} \utens\ofpt{x} & = 0 
    & & \forall \pt{x} \in \Omega\\
    \dTr \utens\ofpt{x} & = \ttens \ofpt{y}  = {g}_N\ofpt{y}
    & & \forall \pt{y} \in \Gamma_N\\
    \Tr \utens \ofpt{x} &= \utens \ofpt{y} = g_D\ofpt{y} 
    & & \forall \pt{y} \in \Gamma_D.
  \end{aligned}
\end{equation}
Here, ${g}_N$ is the prescribed Neumann data in terms of surface tractions and ${g}_D$ represents the prescribed Dirichlet data in terms of displacements.

The variational solution of the BVP~\eqref{eq:BVP} is obtained by means of a boundary integral equation~\cite{hsiao2008}. For linear elasticity, the starting point is
\begin{align}
  \label{eq:somiglianas_identity}
  \Tr \utens\ofpt{x} = \Tr \int_\Gamma \fund{U}(\pt{x},\pt{y}) \ttens \ofpt{y}
  \dgamma{y} - \Tr \int_\Gamma \fund{T}(\pt{x},\pt{y}) \utens \ofpt{y}
  \dgamma{y} & & \pt{y} \in \Gamma.
\end{align}
This equation is known as \emph{Somigliana's identity}, on which the trace operator \eqref{eq:trace} is applied. Hence, all points $\pt{x}$ are shifted to the boundary. In elasticity, $\fund{U}(\pt{x},\pt{y})$ is Kelvin's fundamental solution for displacements and $\fund{T}(\pt{x},\pt{y}) = \dTr_y\fund{U}(\pt{x},\pt{y})$ that for tractions \cite{beer2008}. The subscript $y$ denotes the variable on which the conormal derivative $\dTr$ is applied. The fundamental solutions do not in fact depend directly on the variables $\pt{x}$ and $\pt{y}$ but on their distance $r=|\pt{x}-\pt{y}|$ and are therefore translationally invariant. Moreover, their values behave like $\order(r^{-s})$ with $s\geq 1 \in\N$ in three dimensions and tend to infinity if $r\rightarrow 0$. Hence, the integrals in equation~\eqref{eq:somiglianas_identity} become singular. For the purpose of readability, the boundary integral equation 
\begin{align}
  \label{eq:BIE}
  \operate{(\op{C} + \op{K}) \utens}\ofpt{x} & = \operate{\op{V} \ttens} \ofpt{x}
  && \forall \pt{x} \in \Gamma
\end{align}
is reformulated in terms of boundary integral operators. Here,
\begin{align}
  \label{eq:SL}
  \operate{\op{V} t }\ofpt{x} &= \int_\Gamma \fund{U}(\pt{x},\pt{y}) \ttens \ofpt{y}
  \dgamma{y} & & \forall \pt{x},\pt{y} \in \Gamma
\end{align}
denotes the weakly singular ($s=1$) single layer operator and
\begin{align}
  \label{eq:DL}
  \operate{\op{K} \utens}\ofpt{x} &= \int_{\Gamma} \fund{T}(\pt{x},\pt{y})
  \utens \ofpt{y} \dgamma{y} & & \forall \pt{x},\pt{y} \in \Gamma
  \setminus B_\eps\ofpt{x}
\end{align}
the strongly singular ($s=2$) double layer operator. The latter integral only exists as a \emph{Cauchy principal value}, where the radius $r_{\eps}$ of a sphere $B_\eps$ around $\pt{x}\in\Gamma$ is treated in a limiting process $r_{\eps}\rightarrow 0$. The remainder of that process is an integral free term which is
\begin{align}
  \operate{\op{C} \utens} \ofpt{x} &= c \utens \ofpt{x} && \forall \pt{x} \in \Gamma
\end{align}%
with $c=1/2$ on smooth surfaces. 
\revdel{For non-smooth boundaries, the free term is calculated according to Mantic[39].}{This has been moved to section 3.3. }


%% file: igabem.tex
\section{Discretisation}
\label{sec:IGABEM}

It is generally established that the \emph{isoparametric} philosophy is invoked in isogeometric analysis.
In order to obtain a more efficient BEM formulation, an alternative approach is proposed.
Therefore, the concept of \emph{subparametric} patches is introduced.
They have more control parameters to approximate the variation of the unknowns than for the definition of the geometry and eventually the given boundary data. 
This approach is proposed as a generalisation of the isoparametric concept. 
The independent treatment of the fields leads to a flexible formulation 
which allows a significant reduction of the basis functions for the representation
of given boundary data and the formulation of a collocation scheme which can deal with traction jumps at corners and along edges.

\subsection{B-splines and NURBS}
\label{sec:bspineAndNURBS}
B-splines are piecewise polynomials which are defined by a \emph{knot vector}
\begin{align}
  \label{eq:knot_vector}
  \varXi = \left\{ \uu_1, \uu_2, \dots, \uu_i, \uu_{i+1}, \dots, \uu_n \right\}
\end{align}
which is a non-decreasing sequence of coordinates $\uu_i \leq \uu_{i+1}$ in the parametric space. 
The coordinates \emph{$\uu_i$} themselves are called \emph{knots} and the half-open interval $\left[\uu_i, \uu_{i+1}\right)$ is called the $i$th \emph{knot span}.
The basis functions $N_{i,p}$ are defined recursively, starting with the piecewise constant basis function (order $p=0$), where the support of each $N_{i,0}$ is contained in the $i$th knot span.
\begin{equation}
	    \label{eq:Bspline_N0}
	    N_{i,0}(\uu) = \left\{  
	      \begin{array}{rl}
	      1 & \mbox{ if $\uu_{i}\leqslant \uu < \uu_{i+1}$ }\\
	      0 & \mbox{ otherwise } \\
	      \end{array}\right. 
\end{equation}
Higher order basis functions ($p>0$) are defined by 
\begin{equation}
	\label{eq:Bspline_Np}
	N_{i,p}(\uu)  = \frac{\uu-\uu_{i}}{\uu_{i+p}-\uu_{i}} \: N_{i,p-1}(\uu) 
			+ \frac{\uu_{i+p+1}-\uu}{\uu_{i+p+1}-\uu_{i+1}} \: N_{i+1,p-1}(\uu)
\end{equation}
as a strictly convex combination of basis functions of the previous order $p-1$ \cite{Farin2002b}.
The support of the basis function, $\supp{ \{N_{i,p} \} }=\fromto{\uu_i}{\uu_{i+p+1}}$, is local and entirely defined by $p + 2$ knots.

The refinement of the basis functions $N_{i,p}$ is performed by an extension of the existing knot vector.
The associated refinement procedures are called \emph{knot insertion} and \emph{order elevation}. 
They are comprehensively discussed in \cite{Hughes2005a,Cottrell2009b} and the corresponding algorithms can be found in \cite{Piegl1997b}.
Using knot insertion, a B-spline can be split into \emph{B\'{e}zier segments}.
Thereby, the multiplicity of all knots is equal to the order $p$.
Consequently, the Kronecker delta property is fulfilled at the knots and the basis functions $N_{i,p}$ are the classical $p$th-order Bernstein polynomials which extend over a single non-zero knot span. 

Multivariate basis functions are defined by tensor products of univariate basis functions of parametric direction.
For $d$ dimensions,
\begin{equation}
  \label{eq:BSPLINE_Npq}
  B_{\mIndexI,\mIndexP}(\pt{\uu}) = \prod_{n=1}^{d} N_{i_n,p_n}^n (\uu_n) 
\end{equation}
with multi-indices $\mIndexI=\fromto{i_1}{i_{d}}$ and $\mIndexP=\fromto{p_1}{p_{d}}$. The former defines the position of the basis function in the tensor product structure and the latter represents the order in each parametric direction. 

NURBS are piecewise rational functions, based on B-splines 
which are usually defined by non-uniform knot vectors.
The basis functions are given by 	
\begin{equation}
  \label{eq:NURBS_Rpq}
  R_{\mIndexI,\mIndexP}(\pt{\uu}) = 
  \frac{ 
	  B_{\mIndexI,\mIndexP}(\pt{\uu}) w_{\mIndexI}
	  }{
	  \sum_{\mIndexJ \in K} B_{\mIndexJ,\mIndexP} (\pt{\uu}) w_{\mIndexJ}
	  }
\end{equation}
where $w_\mIndexI$ refers to the associated weight of each basis function and $K$ denotes the index set of all basis functions that are non-zero at $\pt{\uu}$.
Note that $R_{\mIndexI,\mIndexP}$ is equal to $B_{\mIndexI,\mIndexP}$, if all weights are set to the same value.

In this paper, a NURBS \emph{patch} refers to NURBS curves and surfaces which are defined by a particular knot vector.
It is determined by a linear combination of its basis functions
\begin{equation}
  \label{eq:NURBS_mapping}
  \pt{x}\ofpt{\uu} = \sum_{\mIndexI \in K} R_{\mIndexI,\mIndexP}\ofpt{\uu} \: \pt{c}_{\mIndexI}
\end{equation}
where the corresponding coefficients in physical space $\pt{c}_{\mIndexI}$ are referred to as \emph{control points}.
B-spline patches are geometrically located within the \emph{convex hull} of its control points.
To maintain this property for NURBS patches, their weights are chosen to be non-negative.
Since B-splines are a subset of NURBS, the term NURBS patch will be used generally for the remainder of the paper.

\emph{Anchors} are introduced in the following, as a means of linking basis functions to a point at a specific position in the parameter space.
Here, they are generally defined by the Greville abscissa~\cite{Farin2002b} 
\begin{align}
    \label{eq:greville}
    \greville_{\mIndexI} & = \frac{\uu_{\mIndexI+1}+\uu_{\mIndexI+2} + \dots +\uu_{\mIndexI+\mIndexP}}{\mIndexP} .
\end{align}
In the special case of discontinuous (i.e. $C^{-1}$-continuous) basis functions, an additional offset~$\alpha$ is applied to their parametric location.
This offset guarantees that anchors do not coincide and it is defined for a specific $\greville_{\mIndexI}$ by
\begin{equation}
  \label{eq:offset}
  \begin{aligned}
    \alpha_{\mIndexI}  =\frac{
      \sum_{l=1}^{L} \hat{\uu}_{\mIndexI-l}-\hat{\uu}_{\mIndexI} 
      + \sum_{l=1}^{L} \hat{\uu}_{\mIndexI+l}-\hat{\uu}_{\mIndexI}
    }{2L+1} 
    && \with && \hat{\uu} = \vek{\varXi} \bigcup \greville 
    && \und && L = \begin{cases} 1 &\textnormal{if } \left(\mIndexP-1\right) < 2 \\
				  2 & \textnormal{otherwise} \end{cases}
  \end{aligned}
\end{equation}
which is an adaptation of the 2-ring collocation scheme described in \cite{Scott2013a}.
In this $\hat{\uu}$ denotes an extension of the knot vector $\varXi$ in the corresponding intrinsic dimension which includes all Greville abscissae~$\greville$.
The Greville abscissae of discontinuous basis functions are inserted at the position which is closest to its non-zero knot span.
The numbering of the entries depends on the index and position of the current $\greville_{\mIndexI}$, as illustrated in \figref{fig:alpha}.
Note that the Greville abscissae of the discontinuous basis functions $\greville_5$ and $\greville_6$ coincide, which is the initial reason for the application of the offset $\alpha$ for the final anchor.
The anchors for continuous and discontinuous basis functions are depicted in \figref{fig:anchors}.
\tikzfig{tikz/bspline_alpha}{2.5}
{
Definition of the sequence $\hat{\uu}$ for the knot vector $\varXi=\left\{0,0,0,0,1,2,2,2,2,3,3,3,3\right\}$ and the corresponding entries including their indices for the construction of  $\alpha_5$ and $\alpha_6$, respectively.
The thicker ticks denote the positions of the Greville abscissae within $\hat{\uu}$. In the parameter space (top) the Greville abscissae are marked by circles.
}{alpha}{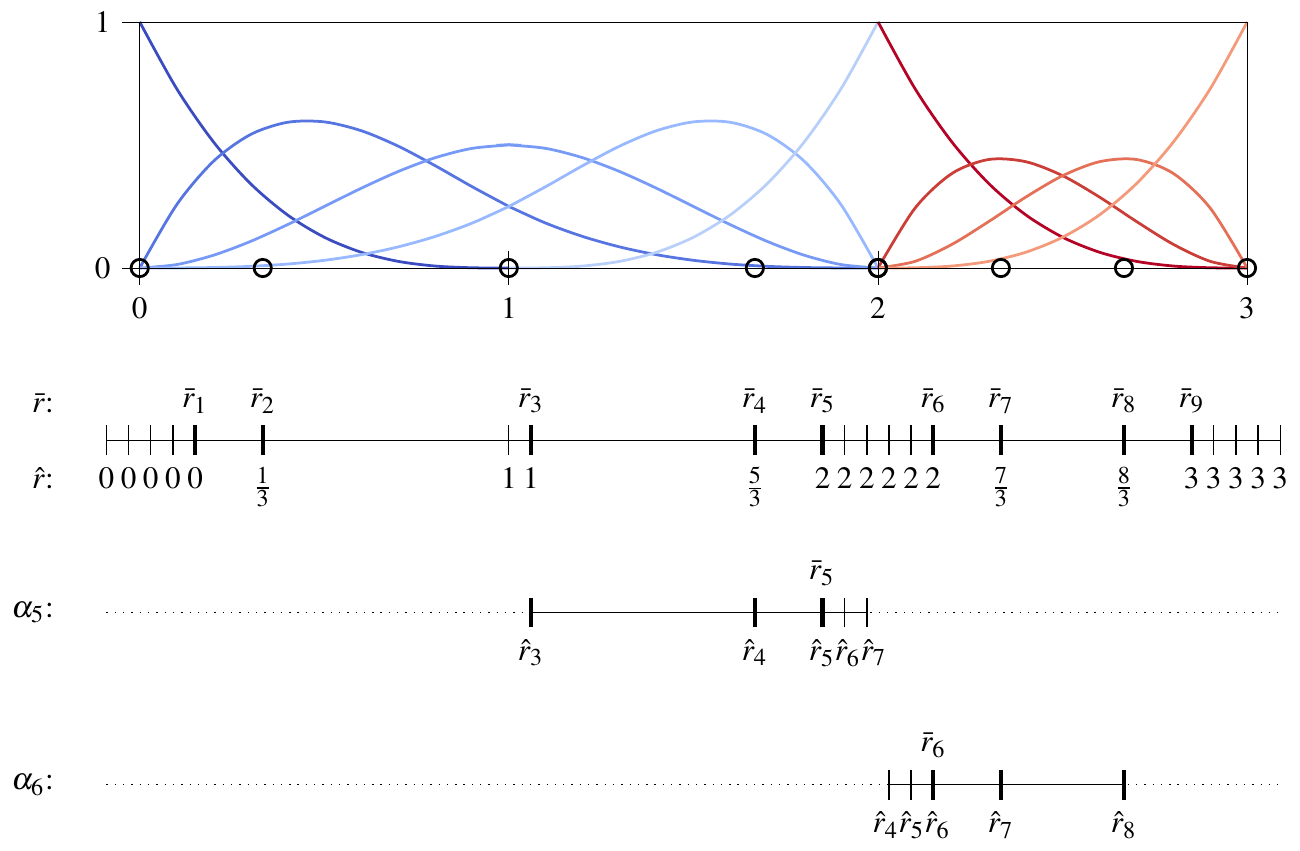}

\tikzfig{tikz/bspline_anchor_new}{1.0}
{
The anchors for the univariate (top), bivariate (bottom), continuous (left) and discontinuous (right) case are marked by black circles.
The arrows indicate the application of an offset $\alpha$. Left: Univariate continuous basis functions (knot vector $\varXi=\left\{0,0,0,1,2,2,3,3,3\right\}$) and 
\revcomment{ 'their corresponding bivariate basis' has been rephrased to point out that the figure illustrates the parameter space due to the tensor product of the univariate basis function}%
the bivariate parameter space constructed by their tensor product.
Right: Univariate discontinuous basis functions (knot vector $\varXi=\left\{0,0,0,1,2,2,2,3,3,3\right\}$) and their corresponding bivariate
\revcomment{ 'basis' has been rephrased to point out that the figure illustrates the parameter space due to the tensor product of the univariate basis function}%
parameter space.
Here, the offset $\alpha$ of the anchors related to discontinuous basis functions is $\pm1/6$ according to \eqref{eq:offset}. 
}{anchors}{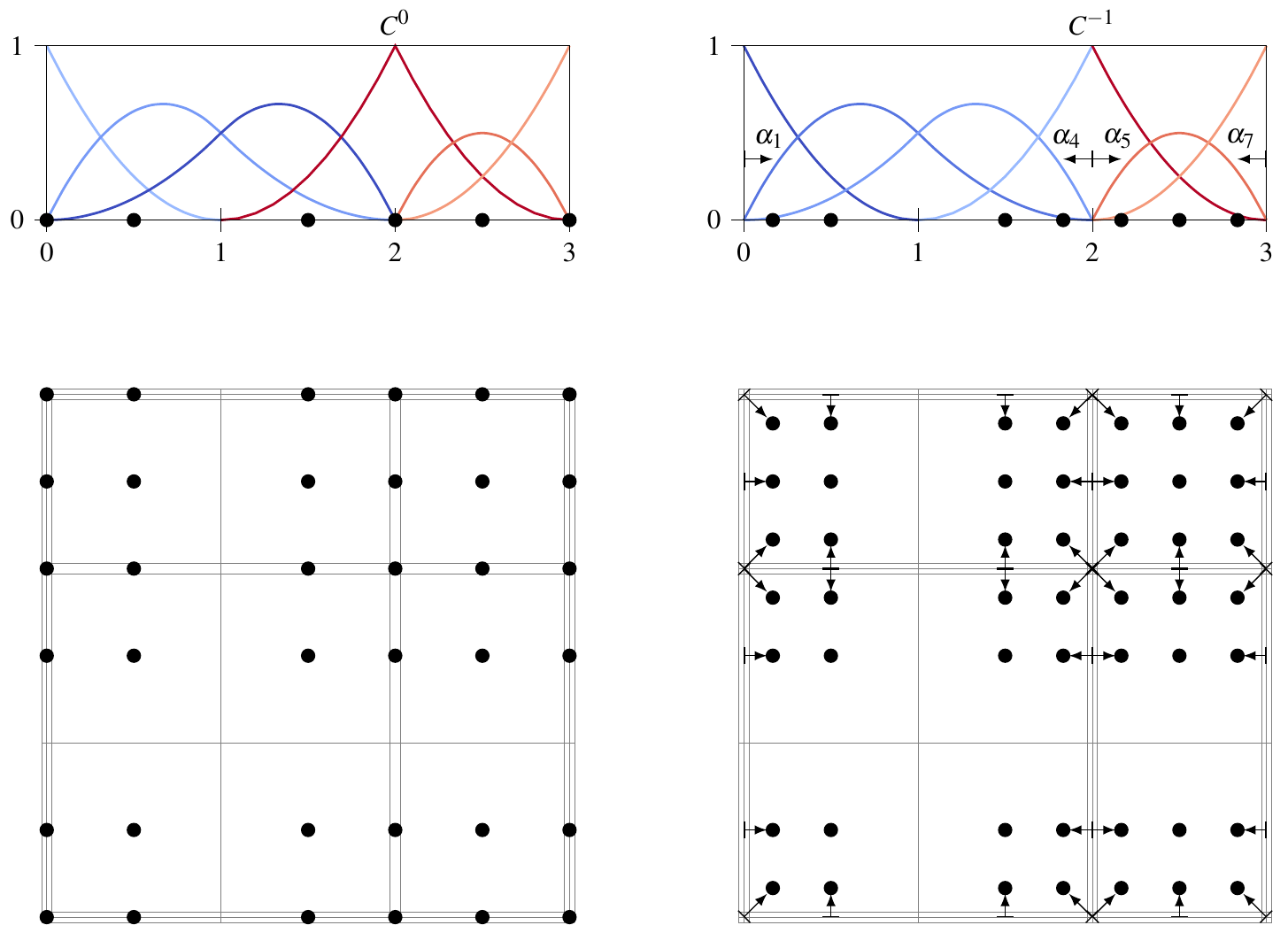}

\subsection{Subparametric Patches}

In this section the concept of \emph{subparametric} patches is introduced. 
It is based on an independent approximation of the geometry, traction and displacement field.
Following the idea of conventional subparametric elements, more parameters are used for the discretisation of the unknown boundary data than for the geometry representation. 

\subsubsection{Geometry Representation}
\label{sec:GeometryRepresentation}

Like most design models, the BEM is based on a boundary representation, hence the parametric dimension is $d-1$ and  
\begin{equation}
  \label{eq:geometry-mapping}
  \mychi\ofpt{\uu}: \R^{d-1} \mapsto \R^d
\end{equation}
is a coordinate transformation, mapping local coordinates $\vek{\uu}=(\uu_1,\dots,\uu_{d-1})^\trans$ of the reference patch to the global coordinates $\pt{x}=(x_1,\dots,x_{d})^\trans$ in the Cartesian system. 
In general, the coordinate transformation is determined by \eqref{eq:NURBS_mapping}.
For the sake of clarity, the multi-index $\mIndexP$ for the order of the basis functions will be skipped and the multi-index $\mIndexI$ unified regardless to the parametric direction in the remaining text.
Thereby the geometric mapping simplifies to 
\begin{align}
  \label{eq:geo_mapping}
  \mychi\ofpt{\uu}:= \pt{x}\ofpt{\uu} = \sum_{k=1}^{K} { \theta_k\ofpt{\uu} } \: \pt{c}_k
\end{align}
where $K$ is now the total number of B-spline or NURBS basis functions $\theta_k$ defined by knot vector  $\varXi_{\theta}$. 

The boundary $\Gamma$ is represented by a disjointed set of NURBS patches $\be$.
\revdel{Hence,}{}%
\revadd{This paper focuses on regular four-sided NURBS patches, hence }{Indicating that this paper deals with patches which have tensor product structure.}%
the computational surface is defined by
\begin{align}
	\Gamma_h = \bigcup_{l=1}^{L} \be_l.
\end{align}
Within each NURBS patch $\be$, any point $\pt{x}$ is evaluated by means of the geometrical mapping \eqref{eq:geo_mapping}. 
Invoking the isogeometric concept, the design model is used directly as computational surface, and thus $\Gamma_h=\Gamma$.
Since this is the most feasible geometry representation available, no geometry approximation error is introduced.
Hence the subscript~$h$ denoting the discretised boundary is dropped in the remaining text.

As a consequence, refinement of $\Gamma$ is unnecessary.
In fact, refinement by knot insertion or order elevation does not change the representation of a NURBS patch, neither in a geometric nor a parametric sense \cite{Piegl1997b}. 
However, the number of elementary operations (i.e. multiplication and division) for the evaluation of geometry data, is directly related to the order of the underlying basis functions.
\revdel{, as illustrated in Figure 3}{This reference has been deleted, because the Figure has been moved to the Appendix.}
A detailed complexity analysis is given in \ref{appsec:appendixOperations}.
It can be seen 
\revdel{in Table 1}{This reference has been deleted, because the Table has been moved to the Appendix.}
that the evaluation is a process with quadratic complexity 
\begin{equation}
  \operations{(p)} = \orderof{p^2}
\end{equation}
related to the order of the basis functions.
Although the result itself does not change, the number of elementary operations increases rapidly with the order.
In particular, for the computation of tangent vectors on NURBS surfaces, which is required for the evaluation of Gram's determinant as well as the calculation of the outward normal.
Due to this inherent increase of computational effort, order elevation of the geometry basis should be avoided.

In a framework with subparametric patches, unnecessary geometry refinement can be completely prevented and the geometry can keep its initial coarse representation, based on the design model, throughout the analysis.
\revcomment{Table 1 has been moved to the Appendix and Fig. 3 has removed from the manuscript.}%

\subsubsection{Discretisation of Cauchy Data}
For elasticity, Cauchy data are vector valued. Hence displacements and tractions are described by
$\uvek=(\utens_1,\dots,\utens_d)^\trans$ and $\tvek=(\ttens_1,\dots,\ttens_d)^\trans$ respectively and discretised by 
\begin{align}
	\label{eq:ansatzDispl}
	\op{Y}_{\utens}\ofpt{\uu} & := \uvek\ofpt{\uu} = \sum_{i=1}^{I} \varphi_i\ofpt{\uu} \: \uvekc_{i}  \\ 
	\label{eq:ansatzTrac}
	\op{Y}_{\ttens}\ofpt{\uu} & := \tvek\ofpt{\uu} = \sum_{j=1}^{J} \psi_j\ofpt{\uu} \: \tvekc_{j}  
\end{align}
inside NURBS patches. The coefficients $\uvekc_{i}$ and $\tvekc_{j}$ are the control parameters of the corresponding field and $\varphi_i$ and $\psi_j$ are the associated basis functions. 
In general, the mappings $\op{Y}$ are similar to the coordinate transformation \eqref{eq:geo_mapping}. 
However, note that the basis functions in the mappings $\op{Y}_{\utens}$, $\op{Y}_{\ttens}$ and $\mychi$ are different from each other.

A unique feature of $\op{Y}_{\ttens}$ is the continuity type of its basis functions $\psi_j$ along non-smooth regions of $\Gamma$.
They are chosen to be discontinuous, in particular $C^{-1}$-continuous, since surface tractions can have jumps at corners and along edges. 
On the other hand, boundary displacements are required to be continuous.
Hence, their basis functions $\varphi_i$ are \mbox{$C^0$-continuous} at non-smooth parts of $\Gamma$. 
This choice of different basis functions for displacements and tractions is a natural choice based on physical constraints.

Without loss of generality, the knot vectors for the Cauchy data are defined, such that
\begin{align}
  \label{eq:knotvectorDisTrac}
 \varXi_{\theta} \subset \varXi_{\varphi} && \und && \varXi_{\theta} \subset \varXi_{\psi},
\end{align}
for the displacement field and traction field respectively.
Hence, the basis functions for the discretisation of the Cauchy data are extended basis functions of the geometry representation. 
Note that the type of the basis functions $\varphi_i$ and $\psi_j$ is not determined by condition \eqref{eq:knotvectorDisTrac}.
For instance, B-spline basis functions can be used for the approximation of the Cauchy data, even though the geometry is represented by NURBS basis functions. 
Such a combination may be more efficient because of the faster evaluation. Additionally, the calculation of new weight values during the refinement process becomes superfluous.

\subsubsection{Individual Refinement}
\label{sec:refinement-strategy}
 
The approximations and related refinements are tailored to the specific needs of the different fields. 
The main idea is to refine only when necessary.
In the proposed approach the geometry field keeps its initial basis which is determined by the design model.
In our implementation, this initial basis is also the origin for the basis functions $\varphi_i$ of the displacement field.
In order to construct the discontinuous traction basis functions $\psi_j$
which exhibits $C^{-1}$ instead of $C^0$-continuities, the initial basis is refined by knot insertion.
An example for an initial basis and its discontinuous version is shown in \figref{fig:initialBasis}.
\tikzfig{tikz/bspline_initial_basis_3d}{1.0}
{
  The initial
  \revcomment{ 'basis' has been replaced by 'parameter space' to point out that the figure illustrates the parameter space due to the tensor product of the shown univariate basis function}%
  parameter space for continuous (top) and discontinuous (bottom) boundary data derived from the knot vector of the geometry representation $\varXi_{\theta}=\left\{0,0,0,1,1,2,2,2\right\}$ in both parametric directions. The corresponding anchors are marked by black circles. 
}{initialBasis}{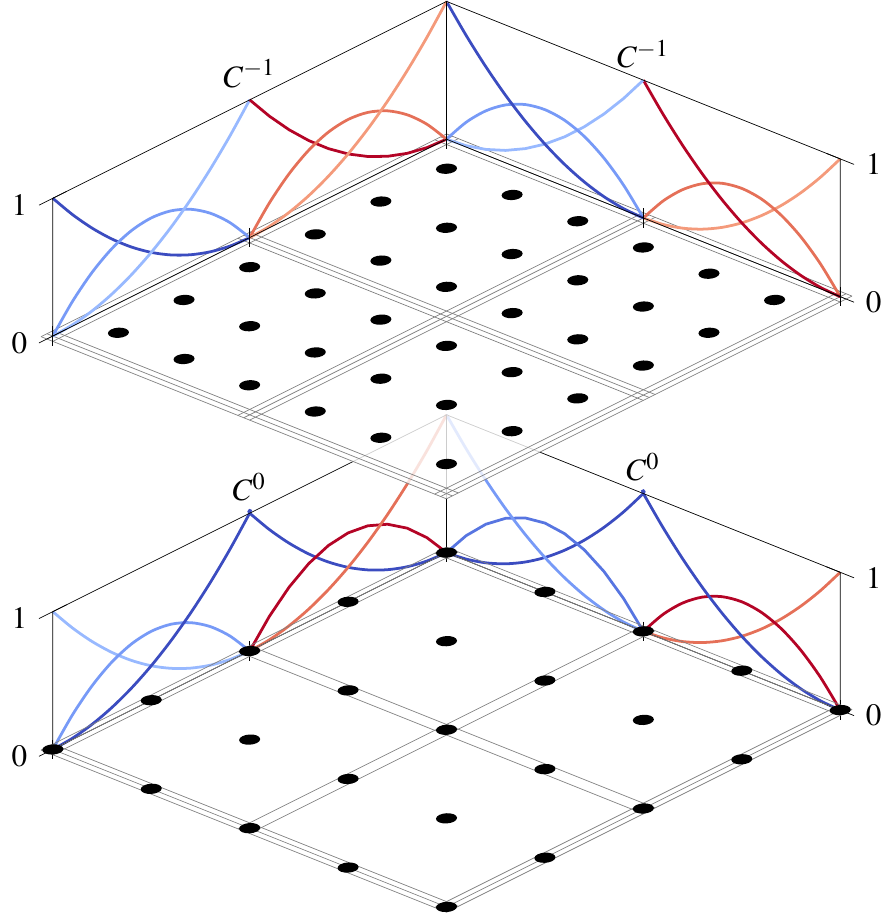}

For the discretisation of known boundary data, two types of boundary condition are distinguished. \emph{Common}~boundary~conditions are exactly represented by the initial basis. Hence, no refinement is required. Examples are fixed displacement and linear traction distribution. Boundary conditions which require a refinement in order to improve their approximation quality are referred to as \emph{complex}.

Usually, the initial basis is not able to describe unknown boundary data sufficiently accurate and therefore, it has to be refined to improve the solution.
Hence, a refinement of the displacement field is performed on the Neumann boundary $\Gamma_N$ and the traction field is refined on the Dirichlet boundary $\Gamma_D$.

To conclude, the geometry keeps its initial representation whereas the basis functions of the unknown fields are refined. 
Hence, there are more basis functions for the unknown fields than for the geometry representation which is the reason for the choice of the term subparametric patch.
In this sense, the proposed approach is a generalisation of the common isoparametric approach used in isogeometric analysis.
For basis functions related to known boundary data the refinement is optional, depending on the type of boundary condition.

\subsection{Collocation}
\label{sec:IGABEM_collocation}
Using a collocation scheme, the boundary integral equation \eqref{eq:BIE} is enforced on a certain set of \emph{collocation points}. 
Each collocation point $\pt{x} \in \Gamma$ generates a discrete equation. 
In order to solve the system of equations, the number of collocation points is determined by the number of unknowns. 

\revdel{
The positions of the collocation points are determined by the anchors of the basis functions related to the unknown field.
}{}%
\revadd{
 Collocation at knots and midpoints of knot spans has been used for uniform basis functions of odd and even orders, respectively~\protect\cite{Arnold1985a,Costabel1992a}.}{}
\revadd{For more general knot vectors, it has been shown that collocation at the Greville abscissae~$\greville$ which are mapped to physical space by $\mychi\ofpt{\uu}$ is an accurate and robust choice for regular NURBS patches with continuous basis functions \protect\cite{Li2011a}.
}{%
We do not explicitly point out that this is for B-spline and NURBS, because it is mentioned before that the term NURBS patch refers to B-spline and NURBS.
}%
\revadd{In order to deal with discontinuous basis functions as well, the anchors of the basis functions defined by~$\greville$ and the offset $\alpha$ are used to determine the positions of the collocation points. This choice guarantees that collocation points do not coincide.
}{}%

\revadd{
In the framework of subparametric patches there are several sets of anchors corresponding to the basis functions~$\theta_k$, $\varphi_i$ and $\psi_j$, respectively.
Regarding collocation the anchors of the basis functions related to the unknown field are essential.
}{}%
Hence the anchors of the displacement approximation \eqref{eq:ansatzDispl} are used on $\Gamma_{N}$ and the anchors of the discontinuous traction approximation \eqref{eq:ansatzTrac} are used on $\Gamma_{D}$.   
\revdel{Using their locations defined by {\protect\eqref{eq:greville}} and {\protect\eqref{eq:offset}} a \emph{mixed collocation} scheme is obtained}{}%
\revadd{This enables a \emph{mixed collocation} scheme}{}
that allows discontinuous tractions along non-smooth boundaries while maintaining the physical constraint of continuous displacements. 

The arrangements of collocation points at corners are shown for different boundary conditions in \figref{fig:collocation}.
Note that the displacement anchor at the join of $\Gamma_N$ and $\Gamma_D$ is not used as a collocation point, because the displacement is known there.

\tikzfig{tikz/mixed_collocation}{0.8}
{
	Mixed collocation scheme at corners of different boundaries based on the anchors of the unknown field. The green colour represents the displacement field whereas the traction field is indicated by purple.
	If an anchor is used as collocation point it is marked by a square, otherwise a circle is used. 	
}{collocation}{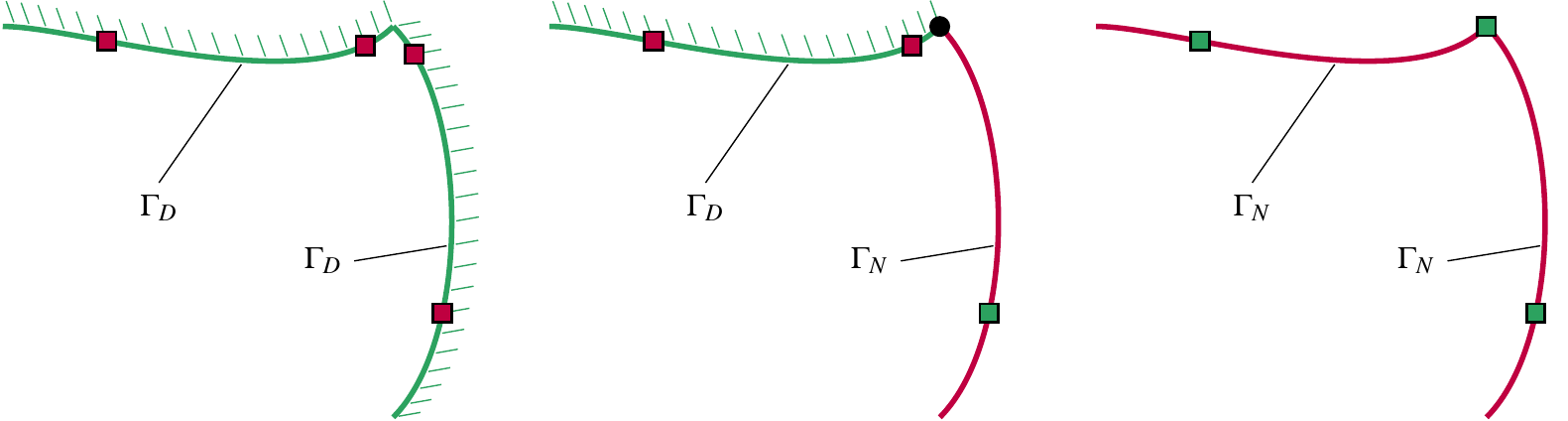}

By partitioning the system of equations with respect to the Neumann and Dirichlet boundary and by shifting the known boundary data to the right hand side, the discretised boundary integral equation can be written in terms of a block system of matrices~\cite{zechner2013}
\begin{align}
	\label{eq:BIEdis}
	\begin{matrix}
		\pt{x}\in \Gamma_{D}: \\
		\pt{x}\in \Gamma_{N}: 
	\end{matrix} \quad
	\begin{pmatrix}
		\vek{V}_{DD} &  -\vek{K}_{DN} \\
		\vek{V}_{ND} & -\vek{K}_{NN}
	\end{pmatrix}
	\begin{pmatrix}
		\vek{\tilde{t}}_{D} \\ 
		\vek{\tilde{u}}_{N}
	\end{pmatrix}
	= 
	\begin{pmatrix}
		\vek{K}_{DD} & -\vek{V}_{DN} \\
		\vek{K}_{ND} & -\vek{V}_{NN}
	\end{pmatrix}
	\begin{pmatrix}
		\tilde{\vek{g}}_D \\ 
		\tilde{\vek{g}}_N
	\end{pmatrix}
	\begin{matrix} 
		\\
		.
	\end{matrix} 
\end{align}
\revadd{It is emphasised that the block-wise setting is essential for the application of \hmatrices{} described in \protect\secref{sec:hmat}.}{Added in order to link the formulation to the application of H-matrices.}

This paper is restricted to NURBS patches which have either Neumann or Dirichlet boundary conditions. Hence, the known control parameters of Dirichlet and Neumann data can be easily evaluated by the inverse of the mappings \eqref{eq:ansatzDispl} and \eqref{eq:ansatzTrac} with
\begin{align}
  \label{eq:knownCauchy}
  \tilde{\vek{g}}_D & = \op{Y}_u^\inv \vek{g}_D && \und &&   \tilde{\vek{g}}_N = \op{Y}_t^\inv \vek{g}_N
\end{align}
respectively. The entries of matrices in the block system \eqref{eq:BIEdis} are
\begin{align}
  \label{eq:BIEdisEntries}
  \vek{V}_{RC}\left[i,j\right] & = (\mathcal{V} \psi_j ) \ofpt{x_i} %
  && \und && \vek{K}_{RC}\left[i,j\right] = ( (\op{C} + \op{K})  \varphi_j ) \ofpt{x_i} 
 \revcomment{Equation has been updated.}%
\end{align}
where $R$ refers to the location of the collocation point $\pt{x} \in \Gamma_{R}$ and $C$ to the location of the integrated basis function $\psi_j,\varphi_j \in \Gamma_{C}$.
\revadd{For non-smooth boundaries, the free term $\op{C}$ at the collocation point is calculated according to \protect\citet{mantic1993}.}{} \revadd{Alternatively, regularisation techniques i.e. as described in \protect\cite{Scott2013a} might be used.}{Information on other approaches calculating the free term has been added.}
\revadd{See~\protect\ref{sec:integration}}{} \revadd{for a more detailed discussion of the numerical evaluation of the matrix entries.}{Referring to the integration part which has been moved to the appendix.}
\revmod{It should be noted}{Minor modification in expression}
that the matrices are fully populated.

In general, the size of \eqref{eq:BIEdis} depends on the number of degrees of freedom $n$ and therefore, the storage requirement of the system matrices is $\order(n^2)$.
Regarding the left hand side of \eqref{eq:BIEdis}, the complexity can only be reduced by the application of a fast solution method, such as the hierarchical matrix approach described in \secref{sec:hmat}.
However, the number of columns~$m$ of the right hand side is determined by the basis functions related to the given data $\tilde{\vek{g}}_D$ and $\tilde{\vek{g}}_N$.
Using subparametric patches, $m$ is independent of $n$ and the refinement described in \secref{sec:refinement-strategy} ensures that it only increases in the presence of complex boundary conditions.
In addition, matrix entries corresponding to homogeneous boundary conditions, such as zero displacements or tractions, are neither calculated nor stored.
In many cases, this leads to a complexity of $\order(nm)$ with $m\ll n$, for the computational effort of the right hand side of \eqref{eq:BIEdis}. 
A comparison of $n$ and $m$ for different boundary conditions is illustrated in \figref{fig:collocRHS}. 

\tikzfig{tikz/collocation_RHS}{1.4}
{
	Comparison of the number of basis functions involved for the representation of the known (above) and unknown (below) boundary data indicated by their anchors (circles) and collocation points (squares).
	The transparent anchors are related to homogeneous boundary conditions and are not taken into account during the analysis.
}{collocRHS}{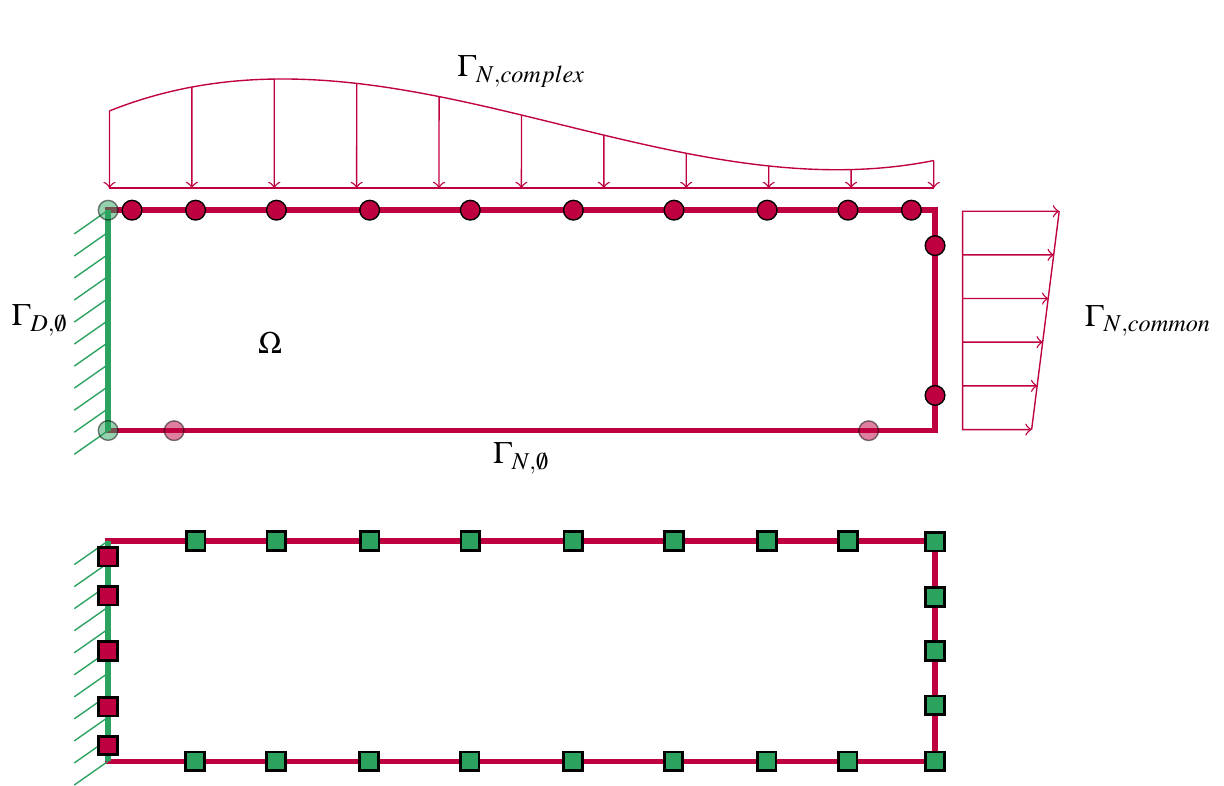}


%% file: hmat.tex
\section{Hierarchical Matrices}
\label{sec:hmat}

As discussed in \secref{sec:IGABEM_collocation}, the system matrices are fully populated and hence, the required storage scales $\order(n^2)$ with the number of degrees of freedom $n$. The same holds for the effort of the involved matrix-vector operations in order to solve the system of equations. In this work, the concept of \hmatrices{} is applied to the isogeometric BEM formulation previously described. 

\hmatrices{} have been introduced by \citet{hackbusch1999} as a general matrix format and applied to the BEM in \cite{boerm2003}. It allows the storage of fully populated matrices with 
\revdel{almost linear complexity}{}%
\revadd{an asymptotic complexity of}{}
 $\order( n \log n )$. Moreover, no further information on the underlying problem is needed after construction of a \hmatrix{}. If the matrix is based on adaptive cross approximation (ACA) introduced by \citet{bebendorf2000}, the underlying physics is also not important for the construction. ACA is a black-box algorithm and hence, has become very popular. Its applicability covers many different problems. 
In the context of the presented formulation, the extension to the elasto-static \cite{bebendorf2006} and elasto-plastic BEM~\cite{zechner2013} is of particular interest.

\subsection{Matrix Approximation }

The reduced complexity of the \hmatrix{} format stems from the low rank approximation of separate blocks $\mat{M}$. It is well known, that the singular values of $\mat{M}$ which is constructed by evaluation of \emph{asymptotically smooth} functions ${\Phi(\pt{x}-\pt{y})}$ decay exponentially if the variables $\pt{x}$ and $\pt{y}$ are well separated. Hence, a degenerated function may be used instead of the original one. To a certain degree of accuracy, the resulting matrix block has numerically low rank \cite{hackbusch2009}. The same holds if matrix entries consist of boundary integrals such as \eqref{eq:BIEdisEntries}. If applicable, $\mat{M}$ may be represented as low rank matrix ($\rk$-matrix) in factorised form
\begin{equation}
\label{eq:low-rank-representation}
\begin{aligned}
  \mat{M}\approx\mat{R}_k = \mat{A} \mat{B}^\trans 
\end{aligned}
\end{equation}
with $\mat{M}\in\R^{m\times n}$, $\mat{A}\in\R^{m\times k}$ and $\mat{B}\in\R^{n\times k}$. If ${k\ll \min(m,n)}$ the storage requirement of $\mat{R}_k$ is considerably reduced compared to the full matrix $\mat{M}$. The number of floating point operations when performing for example matrix-vector multiplication is also reduced.

In this work, ACA is used for the construction of $\rk$-matrices, which is one of the most efficient rank-revealing algorithms available. The construction of $\mat{A}$ and $\mat{B}$ is performed by evaluating chosen rows $\vek{a}_i$ and columns $\vek{b}_j$ of the original matrix block $\mat{M}$. This is the basis for the black-box characteristics of the algorithm. 
Which row- or column-indices are used is controlled by the evaluation of pivot points, typically 
the argument of maximum for the coefficients of a row or column.

\Figref{fig:aca} depicts the construction of a low rank representation by means of the ACA algorithm and shows the reduction of storage of the $\rk$-matrix. Ideally, the number of evaluated rows and columns of the original matrix $\mat{M}$ should be small. How many evaluations are necessary, is controlled by a user-defined approximation quality $\err_{\h}$. As a consequence, the rank $k$ of the approximation is adaptively determined using the criterion 
\begin{align}
  \label{eq:aca_stop}
  \| \vek{a}_{k+1} \|_F \cdot \| \vek{b}_{k+1} \|_F \leq
  \err_{\h} \| \vek{R}_k \|_F && \text{with} &&
  \vek{R}_k = \sum_{\nu=1}^k \vek{a}_\nu \vek{b}_\nu^T
\end{align}
which uses measures in the Frobenius norm $\| \bullet \|_F$. 
The reader is referred to \cite{bebendorf2000} for details of the standard algorithm. The improved ACA+ algorithm used in this work is explained in \cite{grasedyck2005}.%
\tikzfig{tikz/aca}{1.0}{Construction of a low rank matrix $\mat{R}_k$ by ACA}{aca}{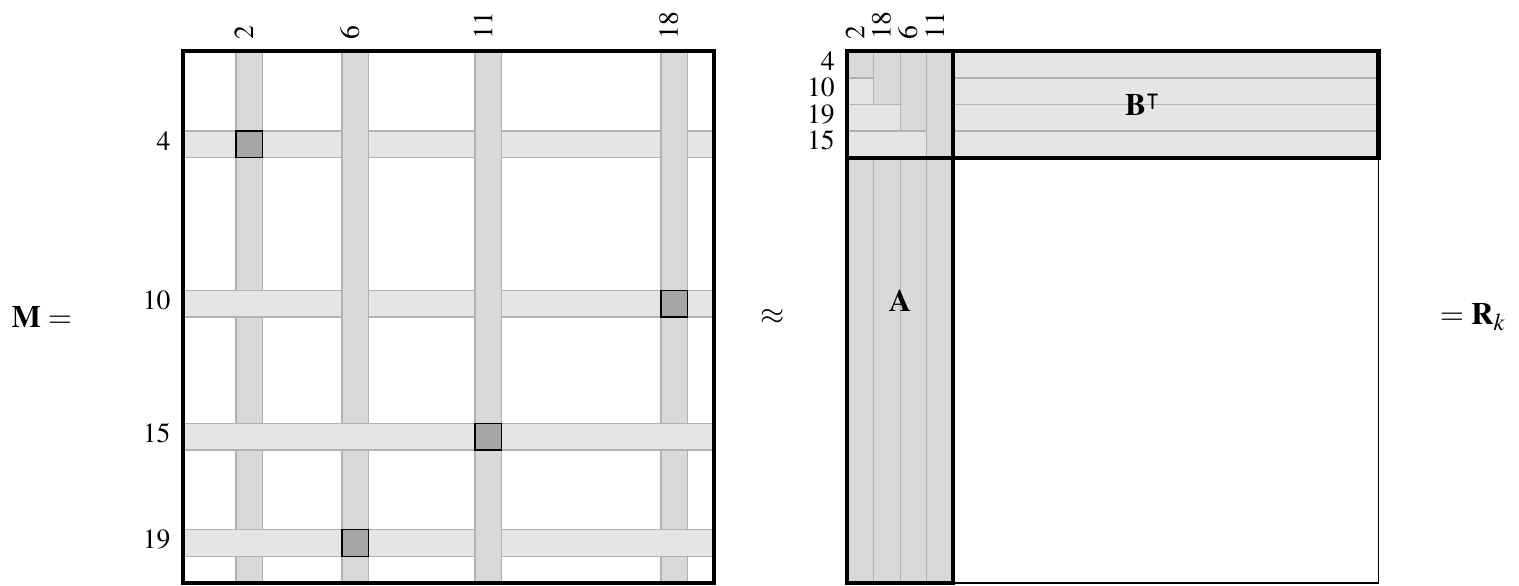}

For problems described by the Lam\'{e}-Navier equation, the fundamental solution is a tensor and therefore vector-valued in each direction. In the context of smoothness, a straight forward assembly of vector valued entries into $\mat{M}$ violates this requirement. Nevertheless, ACA will not fail since the algorithm still reveals the rank of the block, or at least the full rank. But the singular values might not decay exponentially, and thus the rank and the calculation time increase. To overcome this remedy, Cauchy data are decomposed by means of their spatial direction \cite{bebendorf2006,zechner2013} resulting in
\begin{equation}
  \begin{aligned}
    \label{eq:resorting-vectors}
    \tilde{\vek{u}} := &
    \begin{pmatrix}
      \left[ \tilde{\utens}_{1,n} \right]_{n=1}^{N_u} & \dots & \left[ \tilde{\utens}_{d,n} \right]_{n=1}^{N_u}
    \end{pmatrix}^\trans
    \und \\
    \tilde{\vek{t}} := &
    \begin{pmatrix}
      \left[ \tilde{\ttens}_{1,n} \right]_{n=1}^{N_t} & \dots & \left[ \tilde{\ttens}_{d,n} \right]_{n=1}^{N_t}
    \end{pmatrix}^\trans
  \end{aligned}
\end{equation}
for a number of points $N_u$ and $N_t$ in consideration. Consequently, this leads to a block-wise representation
\begin{align}
  \label{eq:block-representation}
  \mat{V} =
  \begin{pmatrix}
    \mat{V}_{11} & \dots & \mat{V}_{1d} \\
    \vdots & \ddots & \vdots \\
    \mat{V}_{d1} & \dots & \mat{V}_{dd}
  \end{pmatrix}\beistrich &&
  \mat{K} =
  \begin{pmatrix}
    \mat{K}_{11} & \dots & \mat{K}_{1d} \\
    \vdots & \ddots & \vdots \\
    \mat{K}_{d1} & \dots & \mat{K}_{dd}
  \end{pmatrix} &&
\end{align}
for the discrete single and double layer operator in $d$ dimensions. Each sub-matrix now contains scalar values based on smooth kernels and ACA can be applied efficiently.

\subsection{Geometrical Bisection}

As mentioned in the previous section, matrix approximation is based on the fact that for asymptotically smooth integral kernels, matrix blocks related to well separated variables $\pt{x}$ and $\pt{y}$ have low rank. A partition of the system matrices with respect to the geometry is thus needed. This means that indices of matrix rows $i\in I$ and columns $j\in J$ are resorted in such a way that their offset corresponds somehow to their geometric distance. The splitting is done block-wise and matrix blocks assigned to be \emph{near field} or \emph{far field}. For the latter type the variables are far away enough from each other and hence, the matrix block is a candidate for approximation.

In the context of \hmatrices{}, this procedure is termed clustering. The geometric properties of row and column indices are defined in two binary cluster trees $\ct$. The procedure is carried out in the following way: first, each index $i$ and $j$ is labelled with a characteristic point. In this formulation this is the collocation point $\pt{x}_i$  or the anchor $\pt{y}_j$ of the NURBS basis function, i.e. $\varphi_j$ as used in \eqref{eq:ansatzDispl}. Second, an axis parallel bounding box $Q$, which includes the support of the corresponding index is generated. For row indices, the box degenerates to a point $Q_i=x_i$, for column indices $Q_j=\supp\{\varphi_j\}$. \revmod{Based upon different strategies, the geometrical bisection which defines the structure of the cluster trees is performed \protect\cite{hackbusch2009,bebendorf2008}}{location of references has been changed}. In the following explanation, the creation of a \emph{geometrically balanced} cluster tree for the row indices $i\in I$ is demonstrated.

Covering all local boxes $Q_i$, a bounding box $B_0^0$ aligned to the initial axis is constructed. In the next step, $B_0^0$ is halved by a cutting plane defined with respect to the direction of its largest extension. The result is now two boxes $B_0^1$ and $B_1^1$, which are then minimised with respect to all $Q_i$ contained within. In the corresponding cluster tree $\ct$, this defines the clusters $t_0^1$ and $t_1^1$ as sons of the root cluster $t_0^0$. The superscript denotes the level $\ell$ in the tree $\ct$ and the subscript is an unique number which relates $B$ and $t$ of the same level. The bisection is continued recursively until
\begin{equation}
  \label{eq:min_leafsize}
  \mathrm{size}( t ) = \# t \leq n_{min}
\end{equation}
is fulfilled, which is characterised by the minimum leaf size $n_{min}$ related to the desired minimal amount of indices in a cluster $t$. The whole procedure is depicted in \figref{fig:clustering} for the minimum leaf size $n_{min}=2$. It can be seen that the diameters of boxes $B$ shrink quickly with increasing level.%
\tikzfig{tikz/clustering}{0.95}{Geometrically balanced clustering for collocation points on a circular geometry and a minimum leaf size $n_{min}=2$.}{clustering}{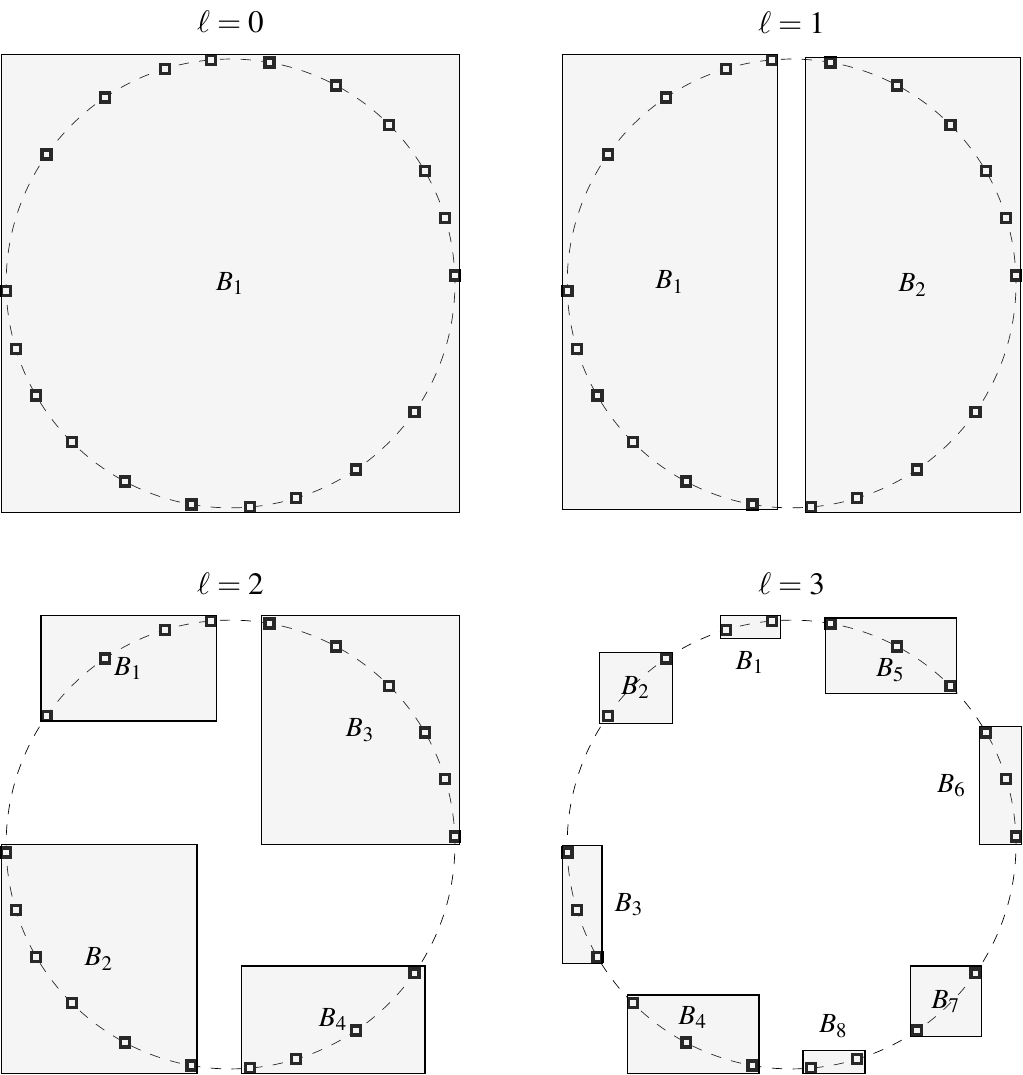}

The same operations are performed for the column cluster tree~$s$, except for the fact that $Q_j$ are now influenced by the support of the assigned NURBS basis function. To provide a reliable bounding box, B\'{e}zier segments are constructed  and the \emph{convex hull} of the resulting control grid is taken. As an example the process is described for a cubic NURBS curve.

The geometry and its control grid is outlined on top of \figref{fig:hull}. The knot vector for the geometry description is defined by $\varXi_\theta=\{ 0,0,0,0,2,4,4,4,4 \}$. Adapted from $\varXi_\theta$, the knot vectors for traction $\varXi_\psi=\{ 0,0,0,0,1,2,3,4,4,4,4 \}$ and displacement field $\varXi_\varphi=\{ 0,0,0,0,0,2,2,4,4,4,4,4 \}$ are defined. Knot vector $\varXi_\psi$ has been subject to knot insertion and $\varXi_\varphi$ to order elevation. The entries of all knot vectors are now accumulated to $\varXi_h$. 
Since geometrical information only is sought, the order of $\varXi_h$ is equal to the order of the geometry description. The B\'{e}zier segments are constructed by means of knot insertions as described in \secref{sec:IGABEM}. 
As a result $\varXi_h = \{ 0,0,0,0,1,1,1,2,2,2,3,3,3,4,4,4,4 \}$.
The corresponding B\'{e}zier segments are depicted on the bottom of \figref{fig:hull}. The resulting control points represent a convex hull of the NURBS curve. Hence, 
a bounding box is generated easily by taking the control points of B\'{e}zier segments related to the support of the considered basis function. For instance, the dashed box in \figref{fig:hull} depicts $Q_1$ for the first basis function $\psi=R_{0,3}$ of the traction description. Whereas $Q_2$ is the bounding box for the support of $\varphi=R_{6,4}$ for the description of displacements.%
\tikzfig{tikz/convexhull}{0.9}{Cubic NURBS curve with its control grid (top) and B\'{e}zier segments (bottom)
  determined by the knot vector for the convex hull $\varXi_h$. The
  dashed boxes $Q_1$ and $Q_2$ denote the bounding box for the support of selected NURBS basis functions describing tractions (red) or displacements (green) respectively.}{hull}{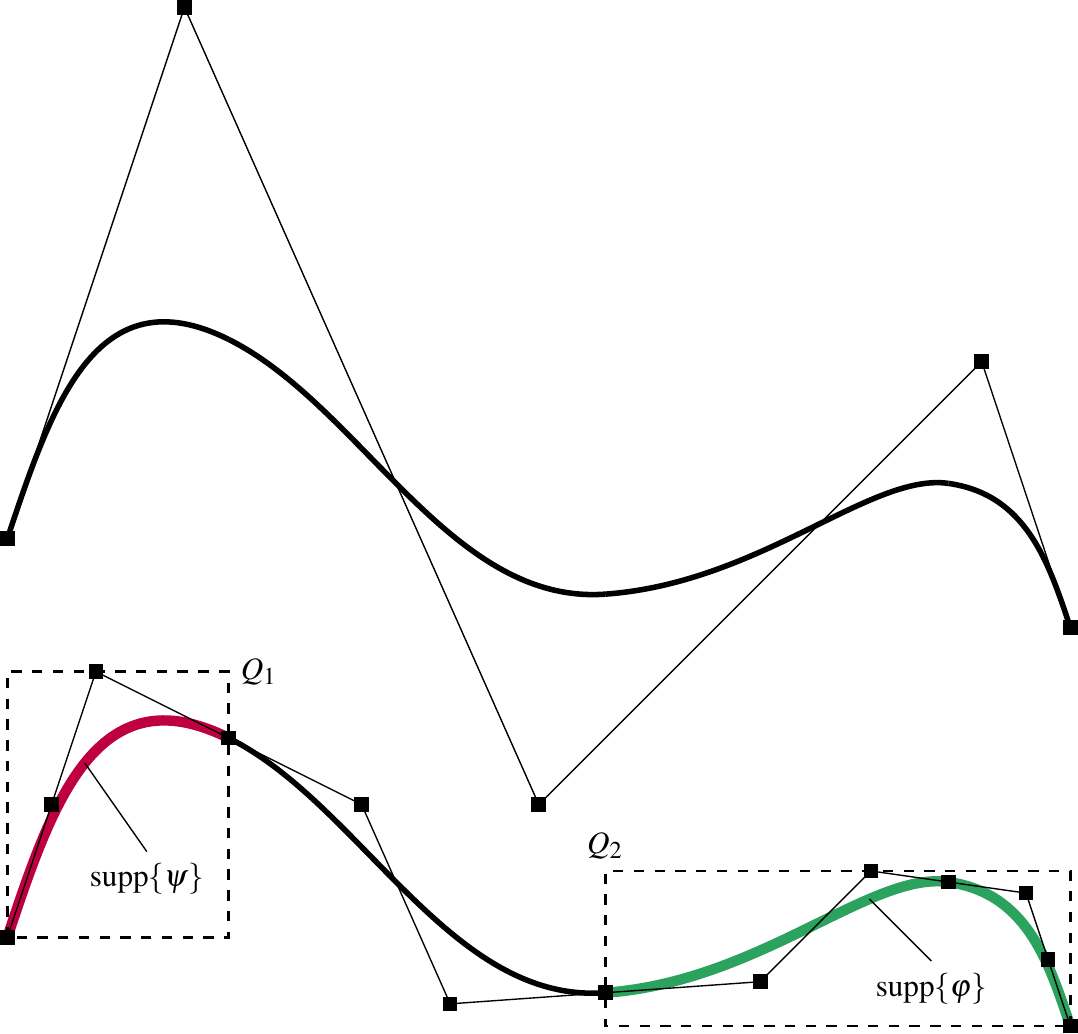}%

The structure of a \hmatrix{} is then defined by the block cluster tree $\bct$ and its nodes ${b = t \times s}$. These nodes contain row- and column-indices taken from the corresponding clusters in the same level $\ell$. The level-wise treatment is not mandatory but reduces the numerical effort. The block cluster tree for off-diagonal matrices in the block system~\eqref{eq:BIEdis} is in fact constructed inhomogeneously by checking admissibility in different levels of the row and column cluster trees. This is because the number of total row or column indices may be relatively unbalanced. An inhomogeneous clustering equilibrates the size of matrix blocks in the \hmatrix{} and thus improves the approximation with ACA.

Efficiency for the evaluation of distances and diameters is the reason for taking axis-aligned bounding boxes for $Q$ and $B$. 
The block cluster tree is a quad tree and forms the basis for the structure of a \hmatrix{}.
An example for the level-wise definition of the matrix structure is depicted in \figref{fig:blockcluster}. 
Starting from $\ell=0$, for each node $b$ an admissibility condition
\begin{equation}
  \label{eq:admissibility-condition}
  \min ( diam(B_t), diam(B_s) ) \leq \eta dist(B_t,B_s)
\end{equation}
is checked \cite{boerm2003}. The variable $\eta$ is called admissibility factor and ranges from $0 < \eta \leq 1$. In effect, the admissibility factor controls the number of $\rk$-matrices: the higher $\eta$ the more admissible blocks are obtained. In \figref{fig:blockcluster} these blocks are coloured green and therefore denote the far field. For red matrix blocks the level in $\bct$ is increased as long as either the row or column cluster is a leaf. Finally, the remaining red blocks not fulfilling \eqref{eq:admissibility-condition} are full matrices and define the near field. These matrix blocks are evaluated with standard BEM techniques whereas far field matrix blocks are $\rk$-matrices to be constructed with ACA.%
\tikzfig{tikz/blockcluster}{1.0}{Matrix partition into blocks defined
  by the block cluster tree $\bct$ up to level $\ell=3$}{blockcluster}{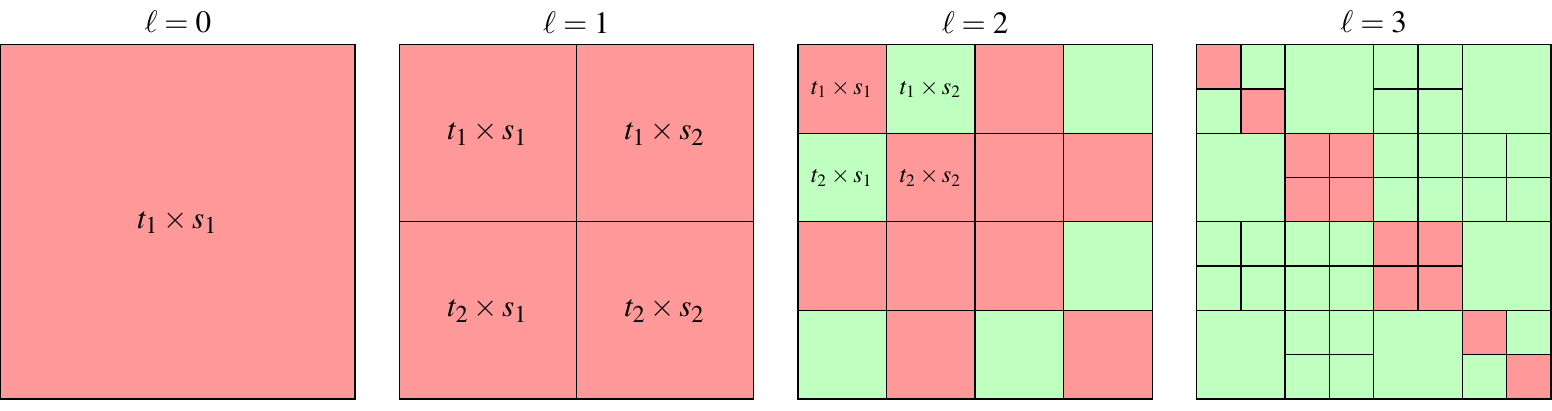}

\subsection{Matrix Operations and Solver}

Compared to other fast methods, \hmatrices{} allow many more matrix operations than merely the matrix-vector multiplication. For example, $LU$ and $QR$ factorization and even inversion are carried out with 
\revdel{almost linear complexity, too}{}%
\revadd{linear complexity up to a logarithmic factor}{} \cite{grasedyck2003}. Of course the accuracy of results is determined by the accuracy of the matrix approximation $\err_{\h}$ itself. For a detailed complexity analysis of operations on and with \hmatrices{}, the reader is referred to \cite{hackbusch2009}.

The complexity of operations is also determined by depth of the underlying block cluster tree. It has been shown by \citet{grasedyck2005}, that after construction the structure of the \hmatrix{} can be coarsened. This is due to the fact that the admissibility parameter $\eta$ could be larger than chosen and that condition \eqref{eq:admissibility-condition} might be sufficient but not necessary. Coarsening yields reduced complexity for matrix operations and less storage requirements as well. In addition, the described approach makes use of the fact that $\rk$-matrices may be re-compressed with only slight additional effort. While a singular value decomposition (SVD) on full matrix blocks is inefficient, it is affordable for the outer product format of $\rk$ matrices. With the help of a truncated SVD, the rank and thus the storage requirement can be reduced further.

The coefficient matrices in \eqref{eq:BIEdis} are ill-conditioned if the number of unknowns becomes very large. 
Thus, special care has to be taken in the solution procedure carried out by \hmatrices{} for large scale computations. As mentioned before, computing the inverse or a $LU$ factorization leads to 
\revdel{almost linear computational complexity}{}%
\revadd{linear complexity up to a logarithmic factor}{}. Hence they can be computed in reasonable time, especially when the approximation error for the particular factorisation is set to a relatively high value, i.e. $\err_\h=10^{-1}$. 
\revdel{For that reason these factorisations can be used as a spectrally equivalent preconditioner for the purpose of solving {\eqref{eq:BIEdis}} iterative with a suitable Krylov subspace method such as GMRES \protect{\citep[see][]{zechner2012b}}.}{}%
\revadd{In the presented implementation, the system of equations is left-preconditioned by such a data-sparse, spectrally equivalent $LU$ factorisation and then solved iteratively with GMRES. For further details on the solution strategy implemented in the described formulation the reader is referred to \protect\cite{zechner2012b}.}{Remark to describe the solution strategy and refer to another paper}


%% file: results.tex
\section{Numerical Results}
\label{sec:results}

In this section, several numerical results for different problems and geometries in two and three dimensions are presented. 
\revadd{The evaluation of matrix entries is performed by numerical integration and an adaptive scheme (see~\protect\ref{sec:integration}), which ensures a user-defined quadrature error~\protect$\err_Q$.}{Introducing the symbol for the integration error}
The superiority of the isogeometric BEM in terms of geometry error compared to standard and higher order Lagrangian boundary elements has already been shown in other publications such as \cite{Simpson2012a} or \cite{Simpson2014a}. Hence emphasis is given to the higher order convergence as well as to the saving of computational costs due to the presented subparametric formulation and the application of \hmatrices{}. 
The presented implementation utilises the \emph{HLib} library~\cite{hlib} for the representation of the \hmatrix{} format and matrix operations.


\input{\CommonPath/results_integration}

\input{\CommonPath/results_convergence}

\input{\CommonPath/results_cantilever3d}

\input{\CommonPath/results_crankshaft}


%% file: results_integration.tex
\subsection{Single Patch Integration}

The following test investigates the efficiency of geometry evaluations with respect to the calculation time. 
Such operations are needed for the calculation of tangents and the outward normal to the NURBS patch in each quadrature point. The test setting is illustrated in \figref{fig:integration_test_setting}. A unit square, defined by a linear NURBS patch $\be$ is located in the $xy$--plane. Three different collocation points $\pt{x}_s$, $\pt{x}_n$ and $\pt{x}_r$ are chosen, so that the boundary integral on $\be$ gets singular, nearly singular and regular respectively. Thus, the integration is subject to different strategies as described in \ref{sec:integration}. The tolerance for the integration error is chosen as $\err_Q=10^{-9}$. The integration is performed for the discrete single layer $\mat{V}_\be$ and double layer potential $\mat{K}_\be$.

\tikzfig{tikz/integration_test_setting}{1.0}{Unit square NURBS patch $\be$ and the location of the collocation points $\pt{x}$ for the regular, nearly singular and singular integration}{integration_test_setting}{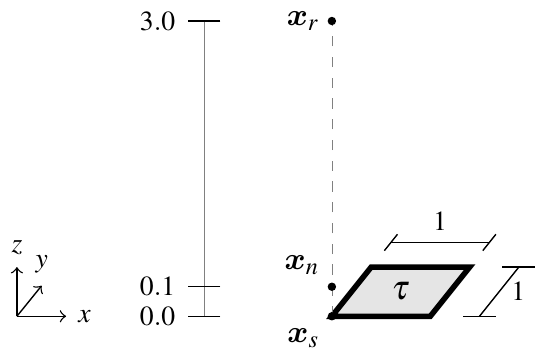}

\revdel{Table 2 lists the number of required Gauss points for this particular case. Due to the subdivision scheme, accurate nearly singular integration takes approximately six times more quadrature points than the regular. While this factor is almost constant for different orders $p_c$ of the Cauchy data basis function, it is envisaged, that $p_c$ plays a significant role when it comes to singular integration. In the defined test setting, precise singular integration demands up to $20$ times more quadrature points than for the regular case. It should be noted that for $\mat{K}_\be$ on a plane, some terms degenerate to zero and accurate integration is thus not as difficult as in the usual case. The number of required Gauss points stays constant as a consequence and is lower than for $\mat{V}_\be$. Nevertheless, accurate integration in the context of a BEM demands a large number of quadrature points and therefore, many geometry evaluations.}{Removed GP table from this section}%
\revcomment{Removed GP table from this section}%

Three different types of patches $\be$ are now investigated. They differ in the description of the geometry and that of Cauchy data. The patch is called \emph{superparametric} if the order of the NURBS basis functions for the Cauchy data is kept constant $p_c=1$ and that of the geometry $p_g$ is refined by means of order elevation. For the \emph{isogeometric} patch, Cauchy data and geometry are described by the same NURBS basis functions. As implemented in the BEM formulation presented, a \emph{subparametric} patch is defined by a geometry representation of lowest possible order, while refinement is performed on the Cauchy data description only. In order to obtain comparable runtimes, the integration process is repeated $100$~times. 

In \figref{fig:time_super} the integration time for different orders on the superparametric patch $t_{sup}$ is related to that of the unrefined patch $t_{0}$. In the case of superparametric $\be$, the number of Gauss points required stays constant, since the NURBS basis functions for the Cauchy data do not vary. \Figref{fig:time_super} thus depicts the increase of computational effort due to more expensive geometrical evaluations with increasing geometry order $p_g$. In the test case and for $p_g=5$ the calculation time is almost twice that much as for the lowest possible order describing $\be$. The evaluation of the traction fundamental solution $\fund{T}$ for $\mat{K}_\be$ additionally requires the evaluation of the unit outward normal $\vek{n}$. Hence the construction of $\mat{K}_\be$ will benefit slightly more, if geometrical information is evaluated on the lowest geometry order available.%

\tikzfig{plots/GP-plane_planestrain_KERNEL_U_and_T_GEOMETRY}{1.0}
{The ratio of integration time on a superparametric and an unrefined patch ($t_{sup}/t_{0}$) for the discrete single layer (left) and double layer (right) operators with regular, nearly singular and singular integration related to the order $p_g$ of geometry description for a fixed number of quadrature points.}{time_super}{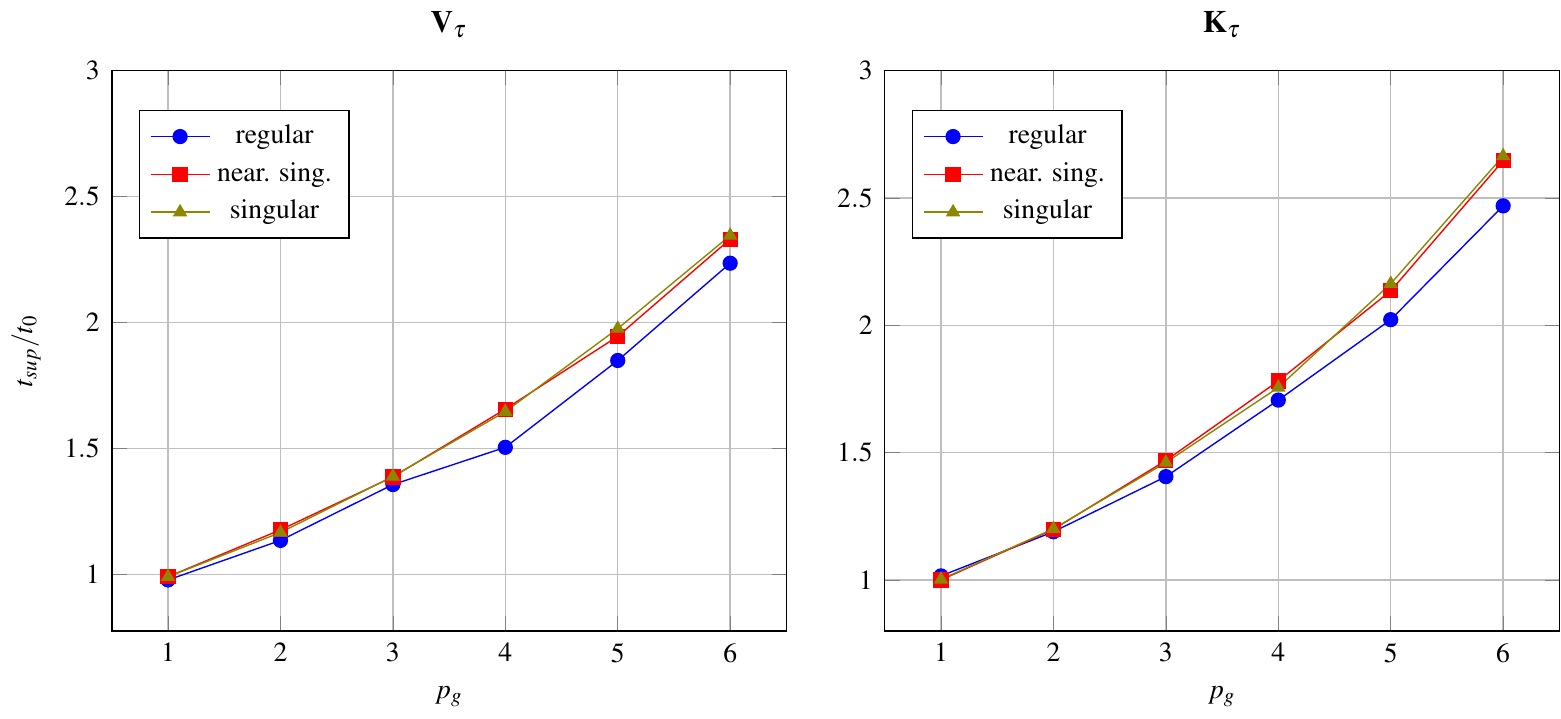}

The presented subparametric formulation makes profit of this in the sense that geometrical evaluations are performed on the original CAD representation. To visualise the performance, the next test case is related to such subparametric patches and hence, the NURBS basis functions for the Cauchy data only are 
\revdel{thus}{}%
refined. In \figref{fig:time_sub} the speedup factor $t_{iso}/t_{sub}$ compared to an isogeometric formulation, where the geometry is also refined, is shown. For the described planar test case, an order elevation by $2$ provides a speedup in patch integration of approximately $25\%$ for $\mat{V}_\be$ and $35\%$ for $\mat{K}_\be$ due to the additional evaluation of $\vek{n}$.%

\tikzfig{plots/GP-plane_planestrain_KERNEL_U_and_T_ALL}{1.0}
{Speedup factor of integration time for the subparametric formulation ($t_{sub}$) compared to the isoparametric BEM ($t_{iso}$) for the discrete single layer (left) and double layer (right) operators with regular, nearly singular and singular integration depending on the order $p_c$ of the Cauchy data.
}{time_sub}{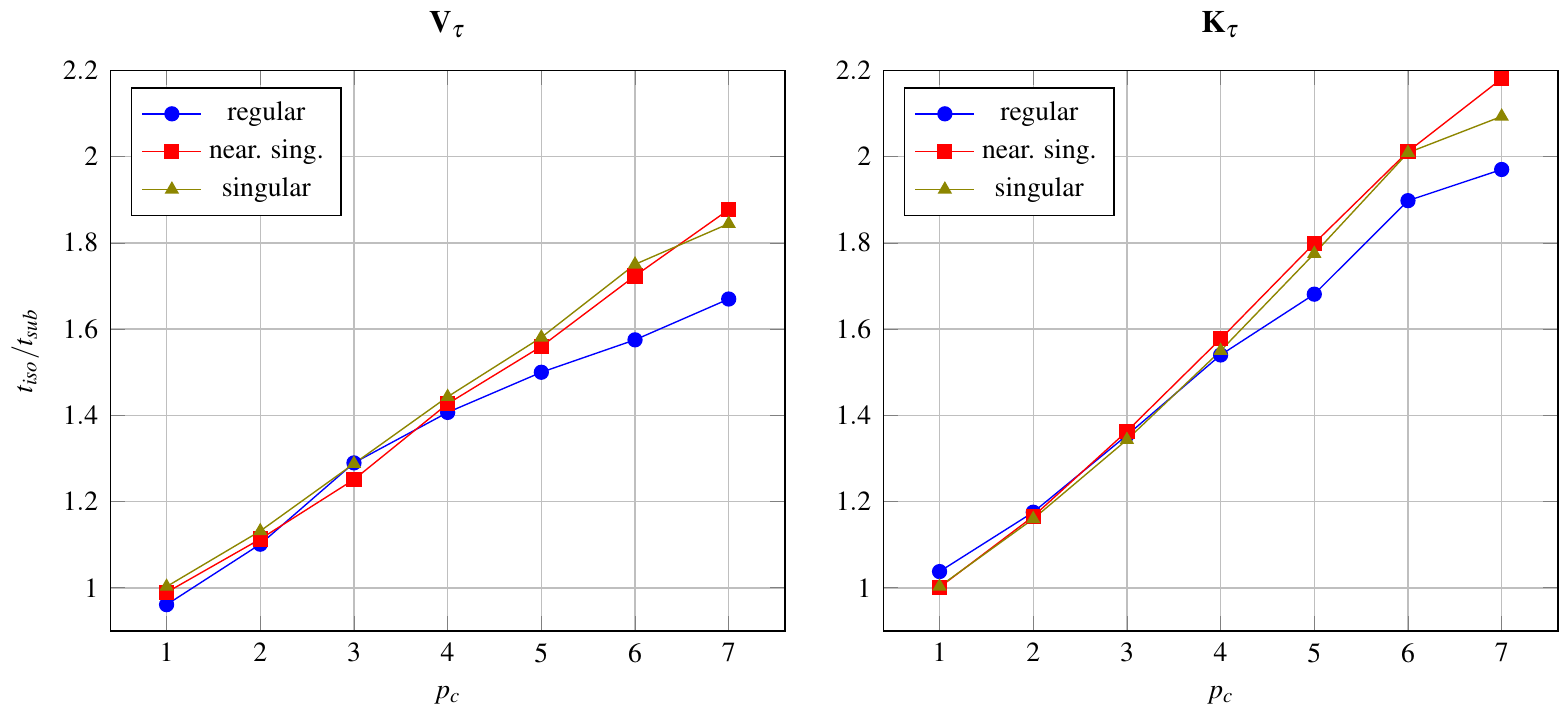}


%% file: results_convergence.tex
\subsection{Convergence Studies}
\label{sec:convergence}

To demonstrate the convergence of the presented implementation on an arbitrary geometry, tests on an infinite domain $\Omega$ with $R\rightarrow\infty$ are performed. The test setting for a circular example is depicted in \figref{fig:test-setting}.%
\tikzfig{tikz/test-setting}{1.0}{Test setting for the exterior problem on a circular geometry with one exemplary source~$\tilde{\pt{x}}$ and field point~$\hat{\pt{x}}$.}{test-setting}{figure13}
To study the accuracy of discretised boundary integral operators, a number of source points $\tilde{\pt{x}} \in \Omega^{-}$ outside the computational domain is chosen, hence $\Omega^{-}\cap\Omega=\emptyset$. From there, the fundamental solution $\fund{U}(\tilde{\pt{x}},\pt{y})$ for a unit load is evaluated at the anchors $\pt{y}$ for the description of displacements on boundary $\Gamma$. 

In order to verify the approximation quality of the discrete integral operators separately, an indirect boundary integral formulation is taken, so that
\begin{equation}
  \label{eq:indirect}
  \begin{aligned}
    (\op{B} \phi) \ofpt{x} & = \utens\ofpt{y} && \with & \utens\ofpt{y} & = \fund{U}(\tilde{\pt{x}},\pt{y}) && \forall \pt{x},\pt{y}\in\Gamma \beistrich \tilde{\pt{x}}\in\Omega^{-} .
  \end{aligned}
\end{equation}
In equation~\eqref{eq:indirect} $\phi\ofpt{y}$ denotes a density on $\Gamma$ which is, in our case, without any relevant physical interpretation. The discrete operator is set to $\op{B} = \op{V}$ for an ansatz with the single layer operator and to $\op{B}=\left( \op{C}+\op{K} \right)$ with the double layer. The density $\phi\ofpt{y}$ is evaluated by solving the system of equations. Consequently, the approximation error $\err_h$ is checked point wise at field points $\hat{\pt{x}}$ inside the domain $\Omega$ with
\begin{equation}
  \label{eq:indirect-check}
  \begin{aligned}
    \err_h & = \utens(\hat{\pt{x}}) - \fund{U}(\tilde{\pt{x}},\hat{\pt{x}})  && \forall\hat{\pt{x}}\in\Omega \beistrich \tilde{\pt{x}}\in\Omega^{-}.
  \end{aligned}
\end{equation}
Interior results are evaluated using the representation formula
\begin{equation}
  \label{eq:indirect-internal}
  \begin{aligned}
    \utens(\hat{\pt{x}}) = (\op{B}\phi) (\hat{\pt{x}}) &&
    \forall\hat{\pt{x}}\in\Omega .
  \end{aligned}
\end{equation}

Furthermore, the overall accuracy of the proposed isogeometric BEM is investigated by the following test setting. The direct formulation 
\begin{equation}
  \label{eq:direct}
  \begin{aligned}
    \left( (\op{C} +\op{K}) \utens \right)\ofpt{x} & = (\op{V} \ttens)
    \ofpt{x} && \with & \ttens\ofpt{y} & =
    \fund{T}(\tilde{\pt{x}},\pt{y}) &&
    \forall\pt{x},\pt{y}\in\Gamma \beistrich
    \tilde{\pt{x}}\in\Omega^{-}
  \end{aligned}
\end{equation}
for the exterior Neumann problem is used. After solving for $\utens$, the approximation error
\begin{equation}
  \label{eq:direct-check}
  \begin{aligned}
    \err_h & = \utens\ofpt{y} - \fund{U}(\tilde{\pt{x}},\pt{y}) &&
    \forall\pt{y}\in\Gamma \beistrich
    \tilde{\pt{x}}\in\Omega^{-}
  \end{aligned}
\end{equation}
 is now evaluated on the boundary only. For 
 an exterior Dirichlet problem, the displacement $\utens\ofpt{y}=\fund{U}(\tilde{\pt{x}},\pt{y})$ is known and \eqref{eq:direct} is solved for $\ttens$ which is checked with the traction fundamental solution $\fund{T}(\tilde{\pt{x}},\pt{y})$. 

These test settings neglect the error of the geometry representation, which for an isogeometric BEM is in any case zero. The relative error is determined by  
\begin{equation}
  \label{eq:relative_error}
  \err_{rel} = \frac{\err_h}{\fund{U}(\tilde{\pt{x}},\pt{y})}
\end{equation}
for the indirect test setting and the direct Neumann problem and by 
\begin{equation}
  \label{eq:relative_error_diri}
  \err_{rel} =\frac{\err_h}{\fund{T}(\tilde{\pt{x}},\pt{y})}
\end{equation}
for the Dirichlet problem.
For the convergence plots, a dimensionless mesh parameter
\begin{equation}
  \label{eq:meshwidth}
  h = \left( \frac{ A^{\ieKs}_{max}}{A} \right)^{1/(d-1)}
\end{equation}
is introduced which plays a similar role than the element length in classical analysis with Finite Elements. The nominator $A^{\ieKs}_{max}$ is the maximal area of all \revmod{non-zero knot spans}{Modified because the 'element' definition is made in the integration section}~$\ieKs$ and the denominator $A$ the complete surface area of $\Gamma$. Equation~\eqref{eq:meshwidth} is used for the tests in two ($d=2$) and three ($d=3$) dimensions.
On smooth surfaces, a rate of convergence for the indirect formulation with the single layer operator of $\orderof{h^{-p-2}}$ may be expected for even orders $p$ \cite{sloan1992}. For odd orders, the double layer ansatz and the direct formulation the rate of convergence should be $\orderof{h^{-p-1}}$ \cite{atkinson1997,sauter2011}. 
\revadd{Usually higher order terms are neglected in error estimates, hence they predict the behaviour for small $h$ only. In general some refinement steps are necessary until the expected rates of convergence appear.}{Added in order to explain the pre-asymptotic behaviour.}


\subsubsection{Tunnel Excavation in Two Dimensions}
\label{sec:convergence_tunnel}

A tunnel excavation which consists of circular arcs is considered as a test case in two dimensions. The geometry and the control points are sketched in \figref{fig:natm2d-geometry}.%
\tikzfig{tikz/natm}{0.8}{Geometry for the tunnel excavation in two dimensions}{natm2d-geometry}{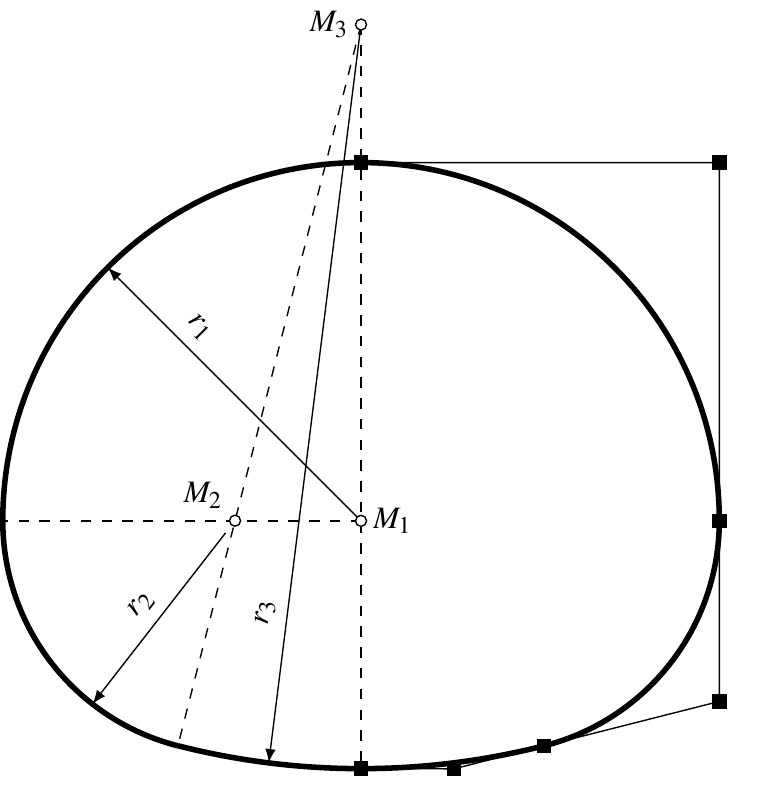}
The measures are $r_1=\SI{4.55}{\meter}$, $r_2=\SI{2.95}{\meter}$ and $r_3=\SI{9.45}{\meter}$. The elastic material properties are $\nu=0.25$ for Poisson's ratio and $E=\SI{10000}{\mega\pascal}$ for Young's modulus. These values are correlated with the Lam\'{e} constants by 
\begin{equation}
  \begin{aligned}
    \lambda=\frac{E\nu}{(1-2\nu)(1+\nu)} & & 
    \text{and} & &
    \mu=\frac{E}{2(1+\nu)} .
  \end{aligned}
\end{equation}
While the geometry description stays constant and is always represented by NURBS basis functions, the discretisation of the Cauchy data is based on either NURBS or B-splines. 
The refinement is carried out gradually by means of uniform knot insertions in the middle of knot spans.
In the studies, which are repeated for different orders ${p=\fromto{2}{5}}$, the appropriate norm of the relative error is plotted against the mesh parameter $h$.

In \figref{fig:natm2d-operators} the accuracy of the discrete integral operators $\mat{V}$ and $\mat{K}$ by means of the indirect formulation \eqref{eq:indirect} is shown.%
\tikzfig{plots/tunnel2d-natm_operators_planestrain_continuous}{1.0}{
  \revcomment{The rate 7 has been corrected to 6}%
  Convergence rates for the discrete single layer (left) and double layer potential (right) on the two dimensional tunnel excavation. 
  \revcomment{Added to clarify the shown gradients.}%
  The optimal convergence rates for the lowest and highest order are indicated by triangles. Note that the rate for the single layer potential is different for even and odd orders. }{natm2d-operators}{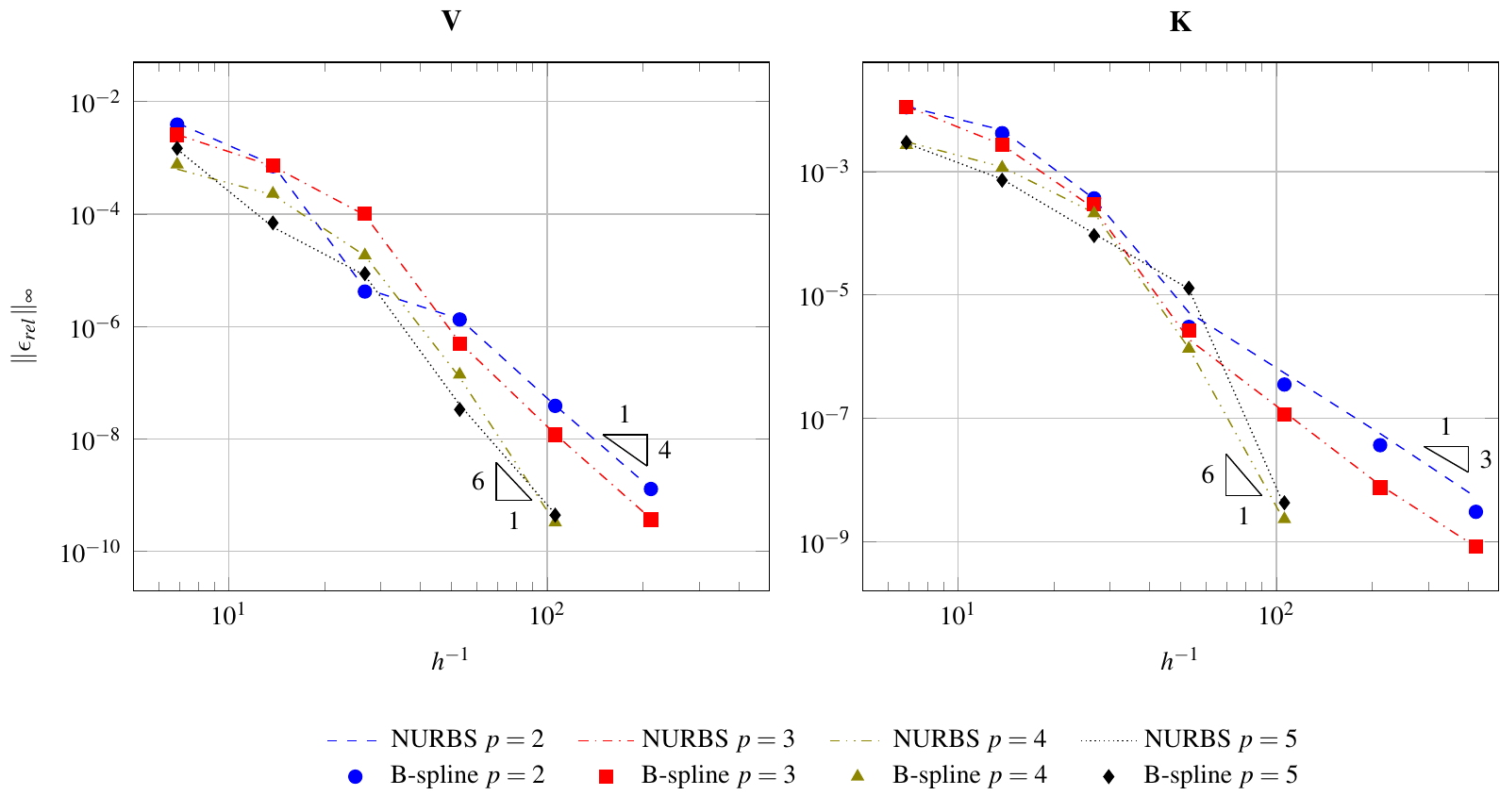}
The error is measured in the maximum norm with respect to the evaluation of $\utens(\hat{\pt{x}})$ in several distributed check points. The integration error is set to $\err_Q=10^{-11}$ and results are shown up to one magnitude less for the relative error. 
Since optimal convergence behaviour is demonstrated for both discrete boundary integral operators, the boundary integral equation and the evaluation of Somigliana's identity for the calculation of internal results are verified. 
In addition, an exterior Neumann as well as an exterior Dirichlet problem is solved with the direct formulation \eqref{eq:direct}.%
\tikzfig{plots/tunnel2d-natm_direct_planestrain_continuous}{1.0}{
  Relative $L^2$-error for the solution of an exterior Neumann (left) and Dirichlet problem (right) on the two dimensional tunnel excavation. 
  \revcomment{Added to clarify the shown gradients.}%
  The optimal convergence rates for the lowest and highest order are indicated by triangles.}{natm2d-direct}{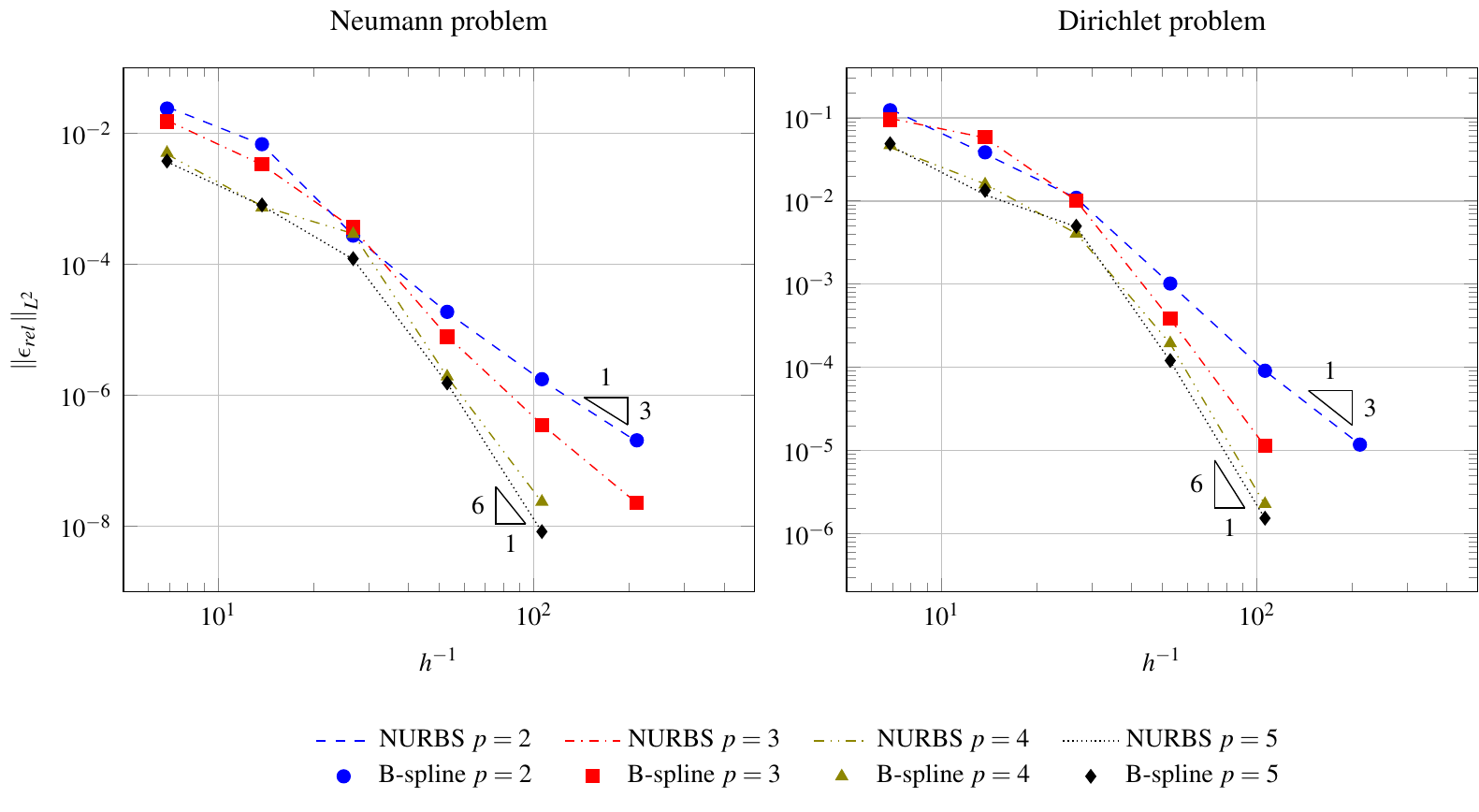}
The relative error is now measured with respect to the $L^2$-norm on $\Gamma$. 
As depicted in \figref{fig:natm2d-operators} and \figref{fig:natm2d-direct}, the approximation of Cauchy data with B-splines and NURBS leads to similar results, as indicated in~\cite{Li2011a}.
Again, optimal rates of convergence are achieved.

The application of B\'{e}zier segments as suggested in \cite{Valle1994a} and \cite{Rivas1996a} is investigated next.
The refinement procedure is thus slightly modified and each new knot is inserted $p$-times.
Exemplary, the results of the direct formulation for the Neumann problem with orders ${p=\left\{2,5\right\}}$ are illustrated in \figref{fig:natm2d-direct-bezier}.
The rate of convergence is similar, but the B\'{e}zier approximation seems to perform better if the relative $L^2$-error is related to the mesh parameter $h$.
The reason for this is, that the number of degrees of freedom within each element is higher than for NURBS or B-splines.
However, if the relative $L^2$-error is related to the total degrees of freedom $n$ the superior accuracy of the smooth NURBS approximation is evident and, independent of the order. This is shown on the right side of \figref{fig:natm2d-direct-bezier}.

\tikzfig{plots/tunnel2d-natm_direct_planestrain_continuous_bezier}{1.0}{Relative $L^2$-error for the solution of an exterior Neumann related to the mesh parameter $h$ (left) and the number of degrees of freedom $n$ (right) on the two dimensional tunnel excavation.}{natm2d-direct-bezier}{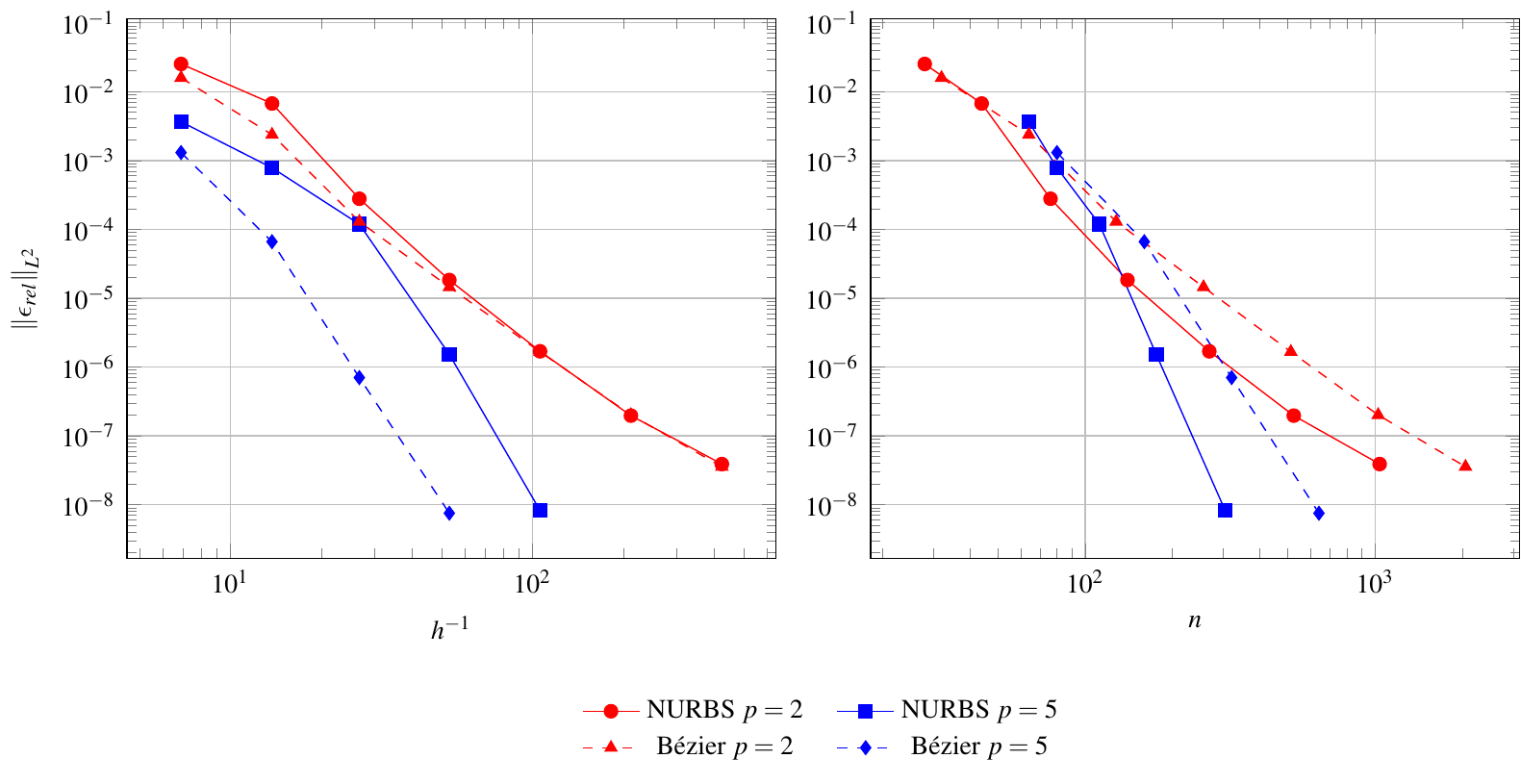}

\subsubsection{Torus}
\label{sec:convergence_torus}

A torus is taken as an example in three dimensions. The geometry is a product of two circles and thus smooth. For the description of its surface, only one NURBS patch is necessary. The shape with major radius $r_m=\SI{5}{\meter}$ and minor radius $r_i=\SI{1}{\meter}$ is outlined in \figref{fig:torus}.%
\tikzfig{tikz/torus}{1.2}{Geometry and measures of a torus}{torus}{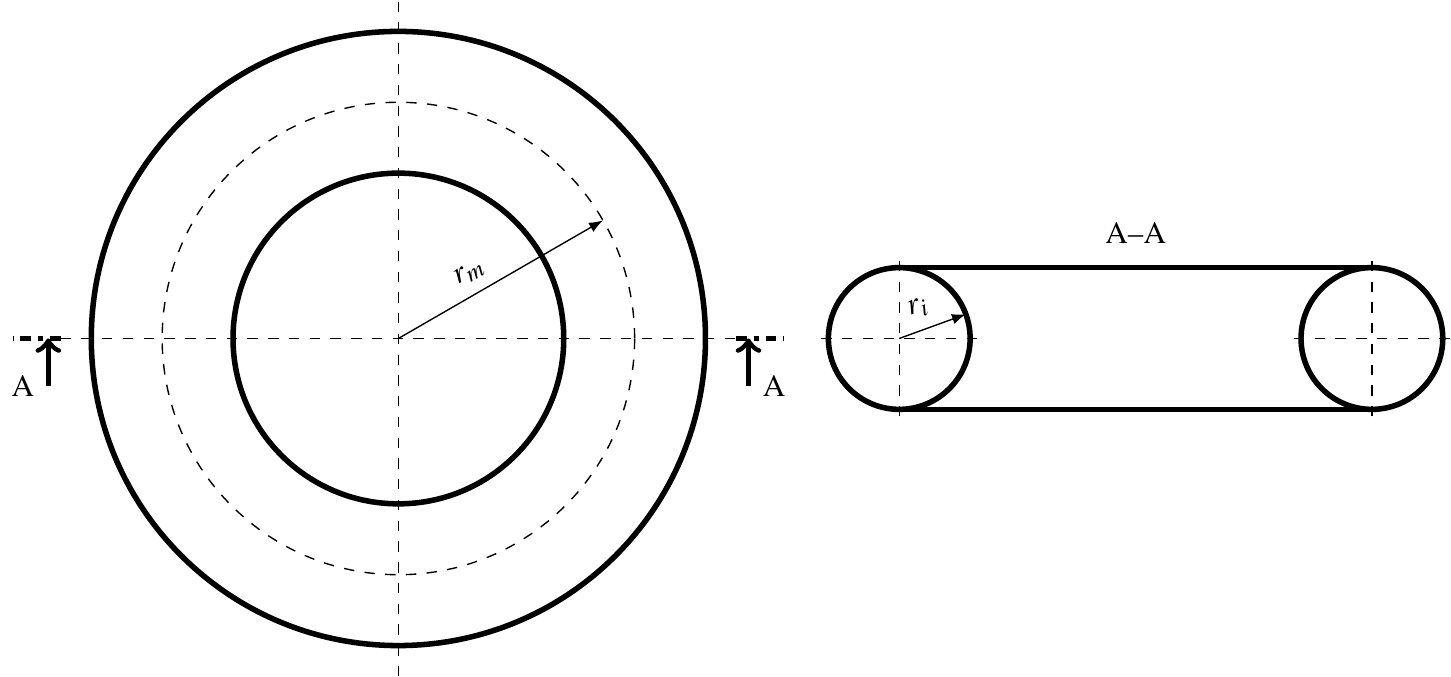}
The material parameters are $E=\SI{1}{\mega\pascal}$ and $\nu=0.3$. 

Convergence analysis for BEM in three dimensions is rather restricted due to the increasing size of degrees of freedom. Moreover, for the Lam\'{e}-Navier equation in three dimensions the displacement $\utens\ofpt{x}$ consists of three components for every spatial point in the domain. Nevertheless, \figref{fig:torus-operators} provides the convergence analysis for the discrete boundary integral operators which is carried out with a system of equations approximated by means of \hmatrices{}.%
\tikzfig{plots/torus_operators_general3d_continuous}{1.0}{
  Convergence rates for the discrete single layer (left) and double layer potential (right) on the torus. 
  \revcomment{Added to clarify the shown gradients.}%
  The optimal convergence rates for the lowest and highest order are indicated by triangles. Note that the rate for the single layer potential is different for even and odd orders. }{torus-operators}{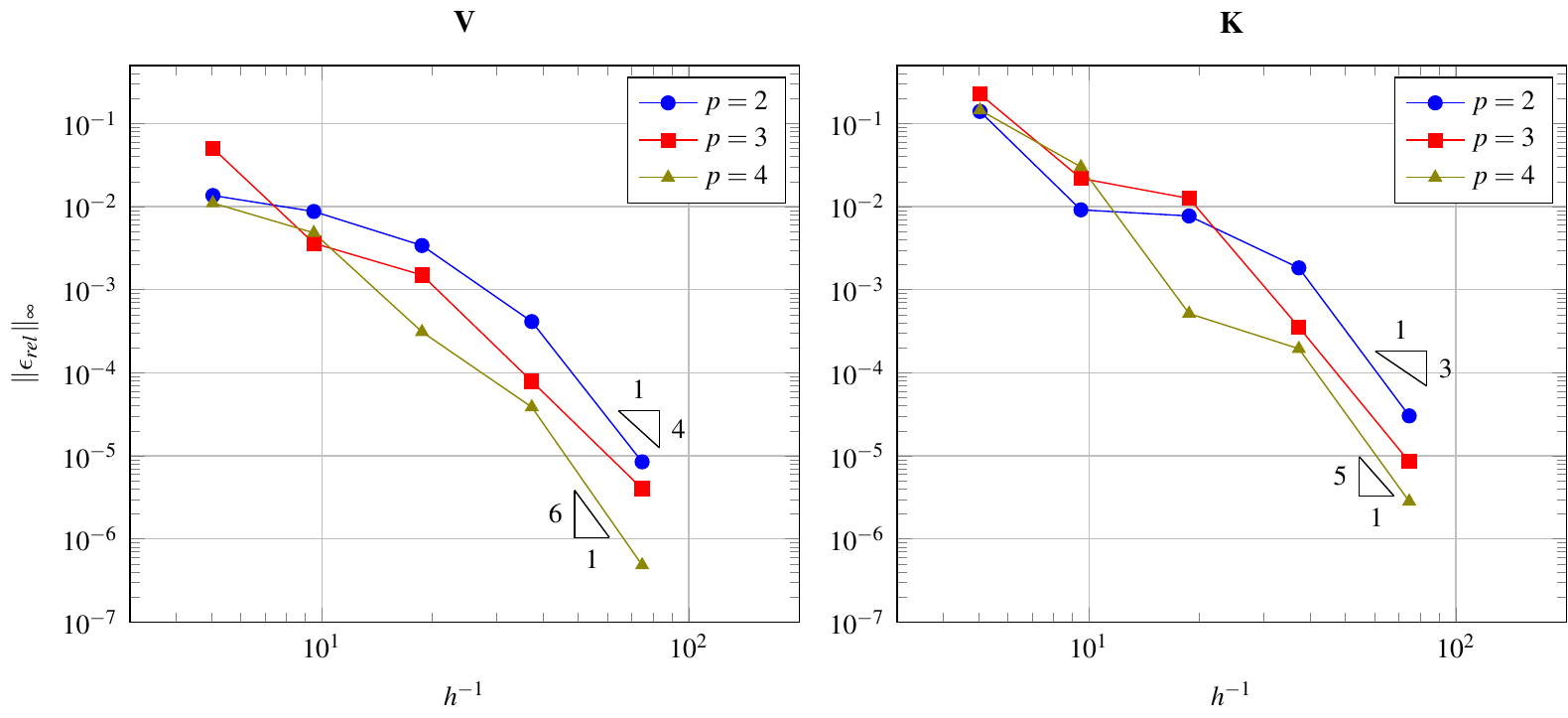}
The accuracy for the ACA algorithm is set to $\err_{\h}=10^{-7}$ and quadrature error one magnitude lower to $\err_Q=10^{-8}$. The admissibility factor is set to $\eta=1$ and the leafsize is adjusted to the size of the problem in order to avoid small matrix blocks. \Figref{fig:torus-operators} indicates, that the matrix approximation has no observable effect on the convergence.%
\tikzfig{plots/torus_direct_general3d_discontinuous}{1.0}{
  Relative $L^2$-error for the solution of an exterior Neumann (left) and Dirichlet problem (right) on the torus. 
  \revcomment{Added to clarify the shown gradients.}%
  The optimal convergence rates for the lowest and highest order are indicated by triangles.}{torus-neumann}{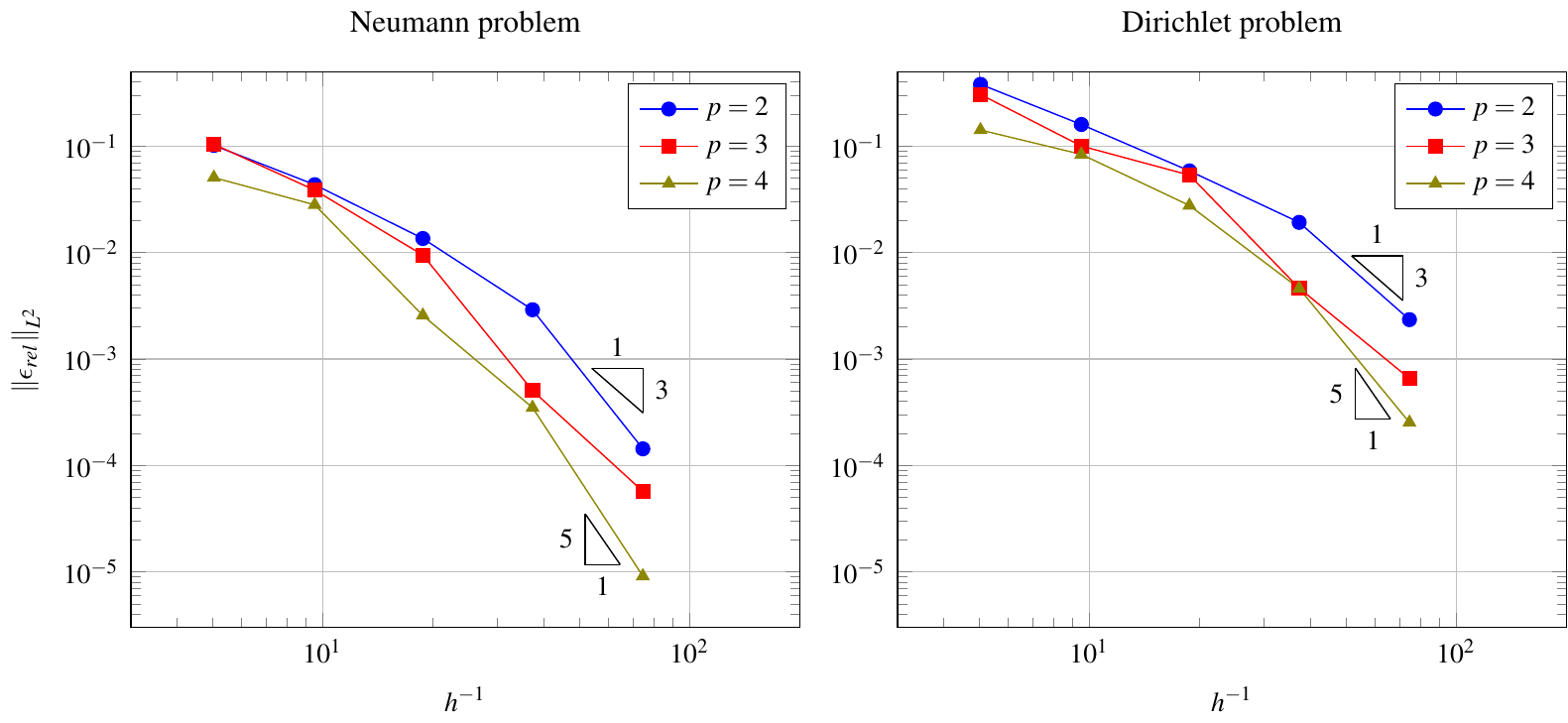}
Finally, \figref{fig:torus-neumann} provides the accuracy for the external Neumann and Dirichlet problem on the same torus measured in the $L^2$-norm. For all tests, optimal convergence rates are obtained.


\subsubsection{Approximation Quality}
\label{sec:convergence_approx}

The implementation presented is especially useful if a certain accuracy is desired. 
\revdel{The error of integration, matrix approximation and that of the residual for the iterative solver are set a priori.}{Removed because the definition of the integration error has been moved to the appendix.}%
As an example, this is demonstrated for the tunnel geometry and the same convergence study as in \secref{sec:convergence_tunnel} for an exterior Neumann problem. But now, the system matrices are subject to \hmatrix{} approximation. The error for ACA is set to $\err_{\h}=10^{-6}$, the error for the numerical integration and for the solver with the same magnitude of $\err_Q=\err_s=10^{-6}$. The admissibility factor and the leafsize are $\eta=1$ and $n_{min}=8$ respectively.

\Figref{fig:natm2d-approx} depicts the relative error $\| \err_{rel} \|_{L^2}$ as well as the compression rate 
\begin{equation}
  \label{eq:compression_rate_h}
  c_{\h}=\frac{\storage( \mat{M} )}{\storage( \mat{M}_{\h} )}
\end{equation}
for the approximation with \hmatrices{} with respect to the total number of unknowns $n$. The factor~$c_{\h}$
relates the storage of the original system matrix $\storage( \mat{M} )$ to that of the approximated  $\storage( \mat{M}_{\h} )$. It is observed, that the relative error for the exterior Neumann problem stays below this threshold for all orders while the system matrices are subject to matrix compression. It should be noted, that the approximation of system matrices by ACA does not influence the convergence until $\err_{\h}$ is reached.%
\tikzfig{plots/tunnel2d-natm_directNeumann_potential_continuous_e5}{1.0}{Convergence and compression for the exterior Neumann problem on the two dimensional tunnel with respect to the degrees of freedom~$n$.}{natm2d-approx}{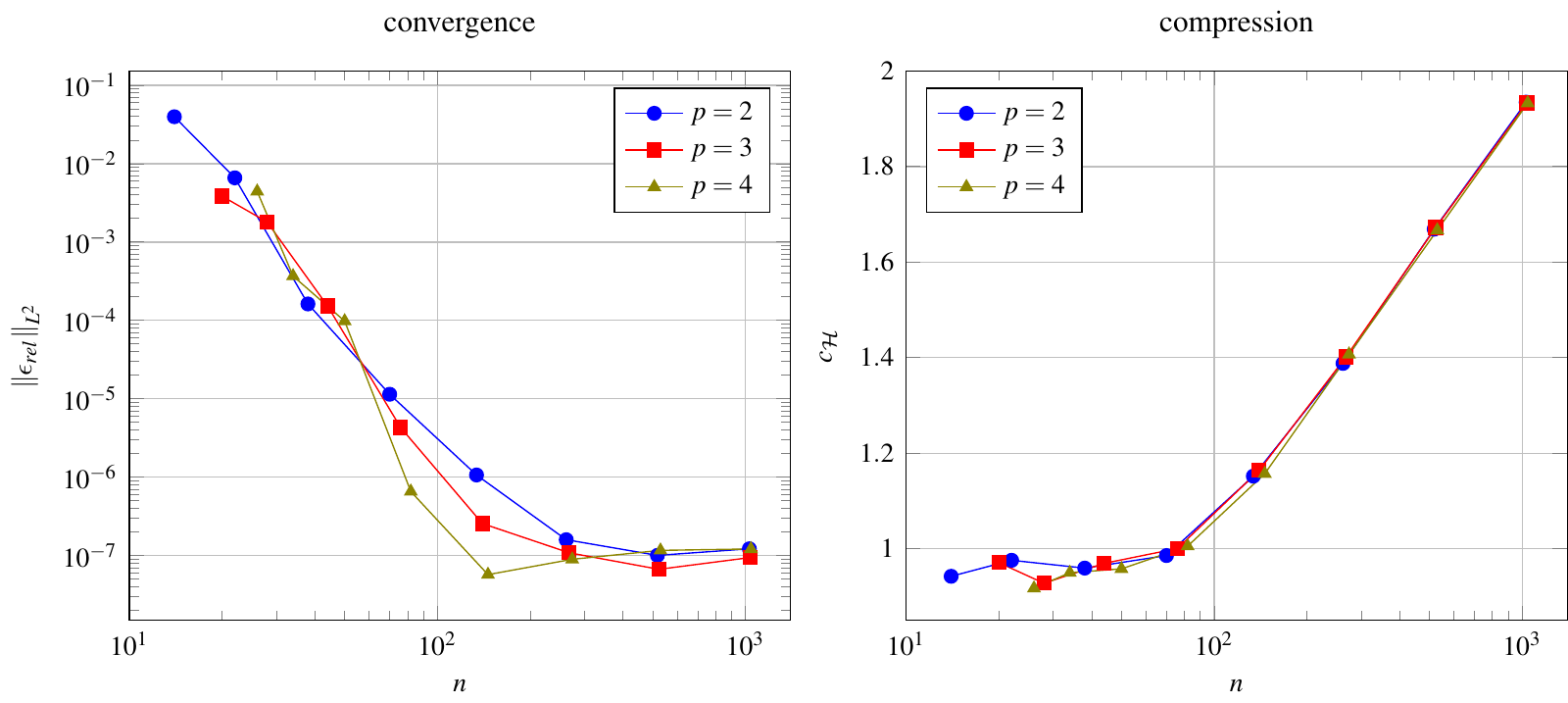}


%% file: results_cantilever3d.tex
\subsection{Cantilever in Three Dimensions}
\label{sec:cantilever}

In order to verify the subparametric approach, a three dimensional cantilever beam with constant vertical load is considered.
The geometry and boundary conditions are shown in \figref{fig:cantilever3D_geometry}. 
The dimensions are $\ell_x=\SI{10}{\meter}$, $\ell_y=\SI{1}{\meter}$ and $\ell_z=\SI{1}{\meter}$ and the elastic material properties are $E=\SI{29000}{\mega\pascal}$ for Young's modulus and $\nu=0.0$ for Poisson's ratio.
At the top, the cantilever is subjected to a constant loading $t_z=\SI{-1}{\mega\pascal}$ in vertical direction. Hence, all other tractions at the free surfaces are zero.%
\tikzfig{tikz/cantilever3D}{1.1}{Geometry and boundary conditions of the cantilever beam}{cantilever3D_geometry}{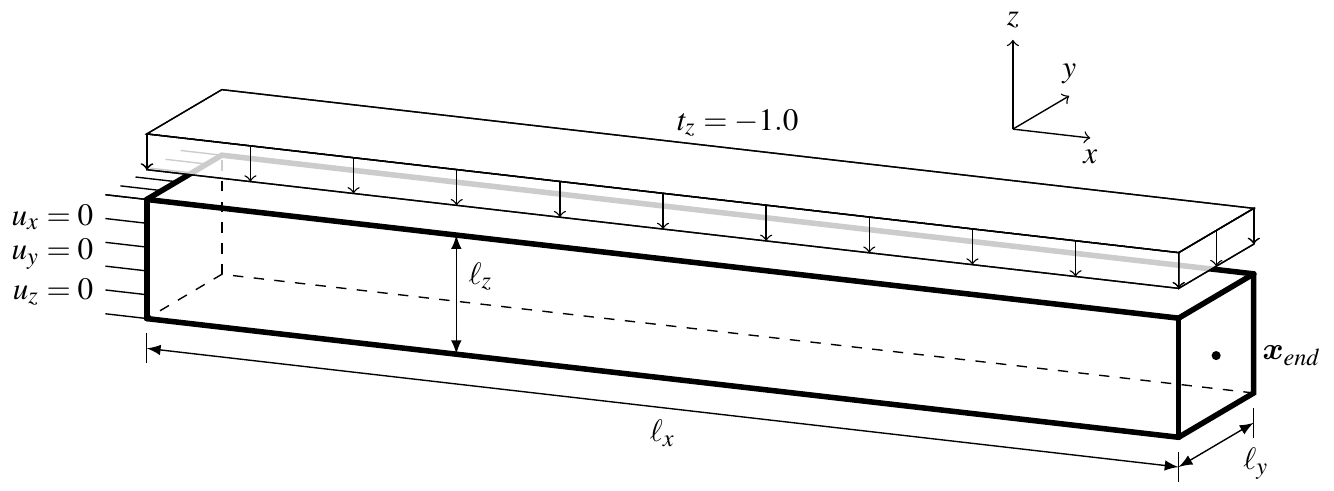}
For different discretisations, the vertical displacement $u_z$ in a point $\pt{x}_{end}$ at the end of the cantilever beam is observed. The test is carried out for both, the subparametric and isoparametric formulation and then compared to the analytic solution by Timoshenko. The results are illustrated in \figref{fig:cantilever3D_displacement}.

Additionally, the compression rate of the right hand side matrix in the block system \eqref{eq:BIEdis}, which is denoted by $\mat{R}$, is investigated. Compression is caused by the subparametric formulation and the avoided refinement with respect to the degrees of freedom $n$ for the known constant tractions. The compression rate relates the storage requirements for the conventional isoparametric approach $\storage (\mat{R}_{iso})$ to that of the subparametric patches $\storage (\mat{R}_{sub})$ and is defined by
\begin{align}
	c_{sub} = \frac{\storage (\mat{R}_{iso})}{\storage (\mat{R}_{sub})}.
\end{align}
Matrix entries related to homogeneous boundary conditions are neglected for the consideration of both, the iso- and subparametric case.%
\tikzfig{plots/cantilever3d_mixed_displacement}{1.0}{
  Left: vertical displacement $u_z$ of the cantilever beam at $\pt{x}_{end}$ for different orders $p$ compared to the Timoshenko beam theory. Right: compression of the right hand side matrix 
  \revcomment{Added to be more specific.}%
  due the application of subparametric patches.
}{cantilever3D_displacement}{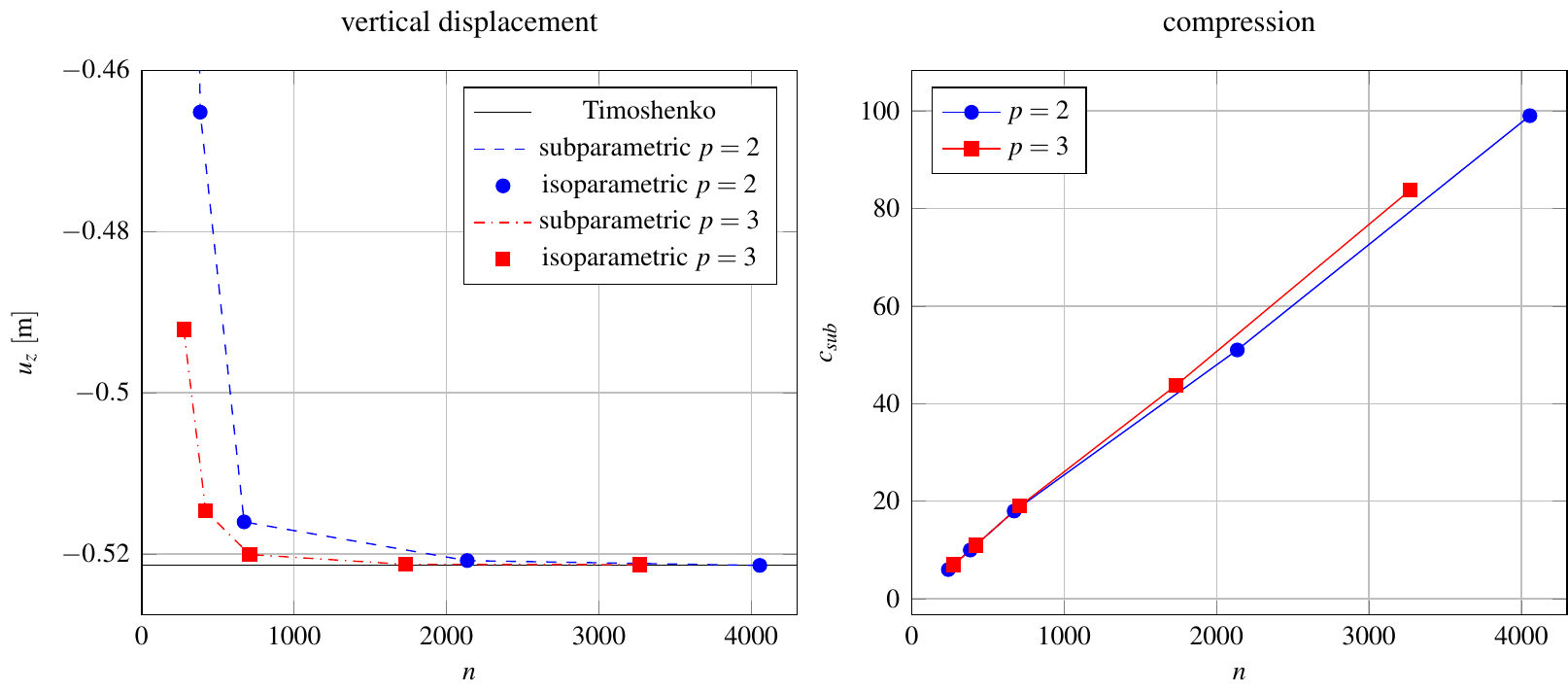}

It is shown in \figref{fig:cantilever3D_displacement}, that the usage of subparametric patches does not alter the quality of results. In terms of storage requirements, the avoided redundancies lead to increasingly high compression rates with respect to the number of unknowns $n$. For the finest discretisation, $\storage( \mat{R}_{sub} )$ is approximately $100$~times smaller than in the isoparametric case.


%% file: results_crankshaft.tex

\subsection{Crank Shaft}

To show the overall performance of the implementation described, a crank shaft is analysed. The geometry is depicted in \figref{fig:crankshaft}.
The material parameters are ${E=\SI{210}{\giga\pascal}}$ and ${\nu=0.25}$. An arbitrary constant load is applied to the crank pins. The flywheel as well as the axle on the other end of the crank are fixed. Because of the constant loading, Neumann boundary conditions on the crank pins are described by means of a coarse discretisation and no refinement is required. All other surfaces are allocated with homogeneous boundary conditions. The calculation is performed with errors ${\err_{\h}=\{10^{-3},10^{-5},10^{-7}\}}$ for the matrix approximation. The minimal leafsize and the admissibility factor are set to $n_{min}=30$ and $\eta=1$.%
\tikzfig{tikz/crankshaft}{1.2}{Crank shaft geometry}{crankshaft}{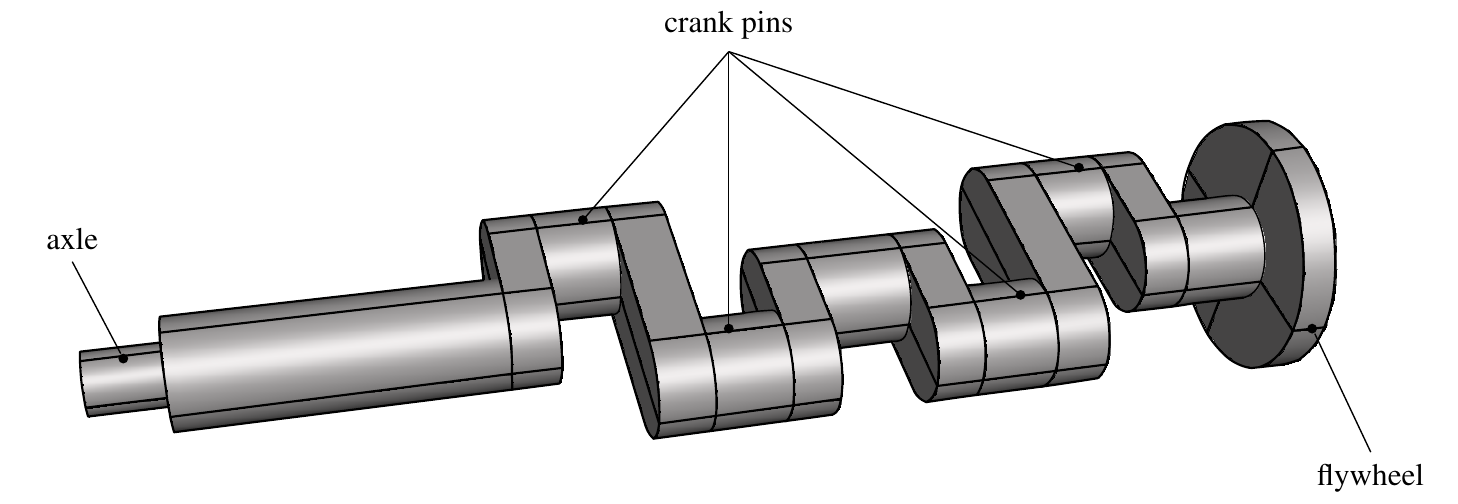}

\Figref{fig:crankshaft_compression} shows the total compression rates $c_{tot}$ for the left hand side $\mat{L}$ and the right hand side $\mat{R}$ of the block system \eqref{eq:BIEdis} with respect to the total number of degrees of freedom $n$. For the matrix $\mat{L}$, which is assigned to the unknown boundary data, the compression is defined by $c_{tot}=c_{\h}$. For the system matrix $\mat{R}$, which is assigned to the unknowns, the total compression rate is calculated by $c_{tot}=c_{\h}\: c_{sub}$. Homogeneous boundary conditions are not considered for the right hand side and the corresponding matrix coefficients are not calculated as for the cantilever in \secref{sec:cantilever}. It is envisaged in \figref{fig:crankshaft_compression}, that for $\mat{L}$ the gradient of the compression rate 
\revdel{becomes almost linear}{}%
\revadd{is of $\order( n \log n )$}{}
with increasing degrees of freedom $n$. Since no further refinement for the description of the known constant tractions is necessary, the number of rows in $\mat{R}$ is determined by the number of collocation points while the number of columns stays constant for all discretisations. The additional compression by means of \hmatrices{} is thus negligible and the slope of $c_{tot}$ is at least linear with respect to $n$.%
\tikzfig{plots/crankshaft_compression}{1.0}{
  \revcomment{The gradient of the figure on the left has been updated to demonstrate the difference compression rates of L and R. }%
  Total compression rates $c_{tot}$ for the system matrices of the crankshaft example with respect to the total degrees of freedom $n$.
  \revcomment{Added to emphasise the difference between the two compressions.}%
  The matrix $\mat{L}$ is solely compressed by the approximation of \hmatrices{} ($c_{tot}=c_{\h}$), whereas $\mat{R}$ is also compressed due to the application of subparametric patches ($c_{tot}=c_{\h}\: c_{sub}$).  }{crankshaft_compression}{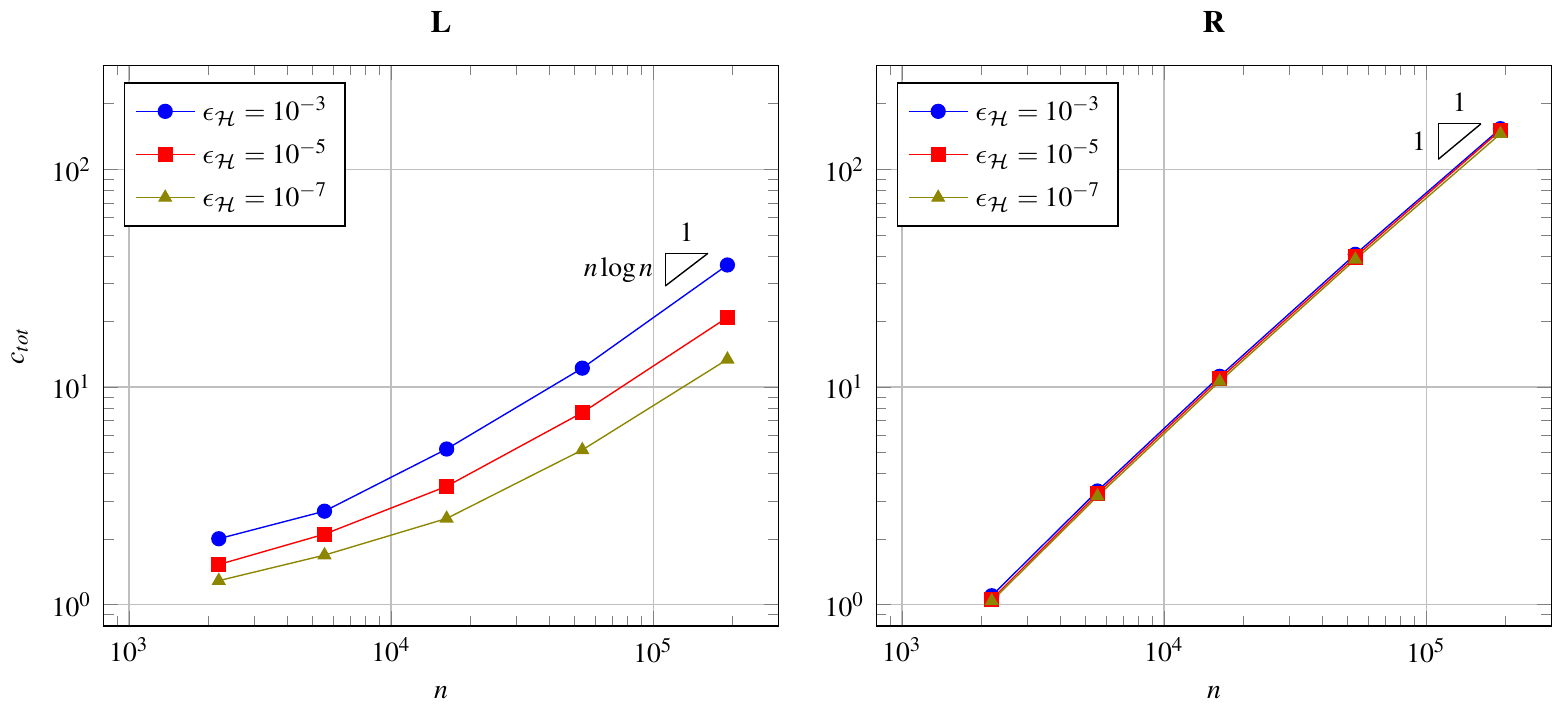}
The resulting absolute displacement for the third level of refinement is illustrated in \figref{fig:crankshaft_displacement}.
\tikzfig{tikz/crankshaft_displ}{1.2}{The absolute displacement $|u|$ of the crank shaft for the third level of refinement.}{crankshaft_displacement}{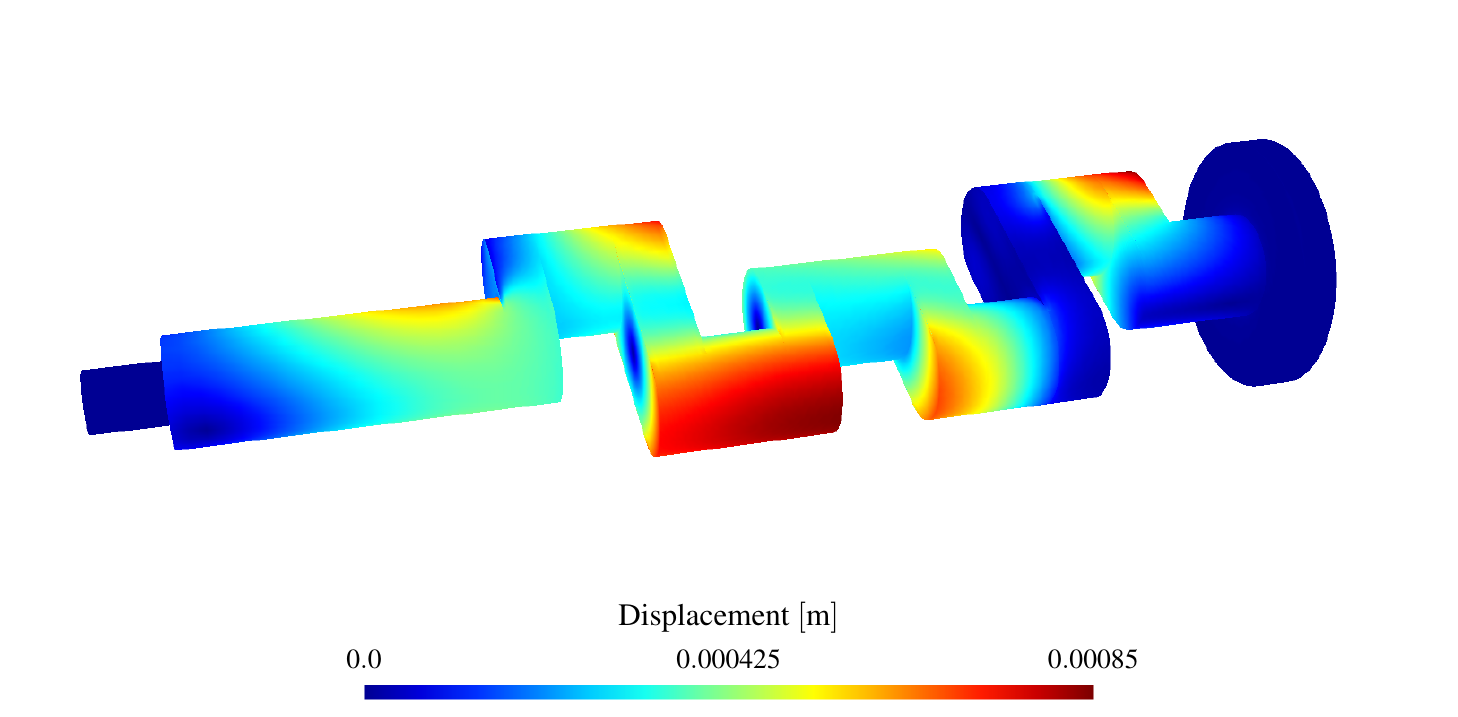}


%% file: conclusion.tex
\section{Conclusion}
\label{sec:conclusion}

A novel isogeometric boundary element method (BEM) based on collocation with NURBS basis functions has been introduced in this paper. The formulation uses an independent description of geometry and Cauchy data and applies the concept of hierarchical matrices (\hmatrices{}). Since the resulting NURBS patches do not strictly follow the isoparametric idea, they are denoted as subparametric patches.

The efficiency of the method presented stems from the fact that the conventional isoparametric concept in isogeometric analysis can be viewed as a process which introduces redundancy: geometry and the Cauchy data are described by more control parameters than actually necessary.
This superfluous information is solely introduced because refinement is mandatory for a sufficient accurate approximation of the solution. 
Refinement is performed independently for each field with subparametric patches and redundant information is thus avoided.
The separation of the field discretisations is beneficial, especially in the context of isogeometric analysis, since known boundary data are defined over a whole patch rather than as an element in standard BEM formulations.
This results in a significant reduction of computational effort for the evaluation of the right hand side of the system of equations.
Additionally, the quadratic computational complexity of the left hand side is reduced to 
\revdel{almost linear complexity}{}%
\revadd{an asymptotic complexity of $\order( n \log n )$}{}
by means of \hmatrices{}. This is achieved by a block-wise setting of integral equations which leads to strictly separated discrete integral operators. The system matrices are approximated by means of adaptive cross approximation. The necessary geometrical bisection uses the convex hull property of NURBS patches with respect to their control points. 

Several numerical tests in two and three dimensions confirm the optimal convergence and accuracy of the presented method, as well as the potential reduction of storage requirements and computational effort compared to conventional isogeometric BEM formulations.
The approximation of the Cauchy data with B-splines or NURBS leads to almost identical results.
Using non-smooth B\'{e}zier segments reduces the accuracy with respect to the number of degrees of freedom. This approximation is not to be mistaken with B\'{e}zier extraction procedure.

\revdel{An adaptive integration scheme has been developed to ensure the accuracy of the integration. 
In order to improve the efficiency of the numerical integration more sophisticated strategies should be developed.}{Removed because the integration section has been moved to the appendix.}%
The choice of discontinuous collocation points presented here is suitable for the examples shown.
However, the issue of their optimal position needs to be investigated in further detail.

\section{Acknowledgment}

This research was supported by the Austrian science fund FWF, provided under Grant Number P24974-N30. This support is gratefully acknowledged.


%% file: appendixDerivatives.tex
\section{Basis Function Derivatives}
\label{appsec:appendixDerivatives}

The first derivative of the B-spline basis functions is computed by
\begin{equation}
	\label{eq:Bspline_Np_Deriv}
	N^\prime_{i,p}(\uu)  = \frac{p}{\uu_{i+p}-\uu_{i}} \: N_{i,p-1}(\uu) 
				    - \frac{p}{\uu_{i+p+1}-\uu_{i+1}} \: N_{i+1,p-1}(\uu)
\end{equation}
and for NURBS basis functions\footnote[1]{The original article contains an error in the definition of NURBS derivatives. This has been corrected in equation \eqref{eq:NURBS_Rp_Deriv} and \eqref{eq:NURBS_Rpq_Deriv}.} it is given by
\begin{align}
	\label{eq:NURBS_Rp_Deriv}
	 R^\prime_{i,p}(\uu)&=\frac{N^\prime_{i,p}(\uu) w_{i} - 
	 N_{i,p}(\uu) w_{i} \beta
	  }{\sum_{j=0}^{n}N_{j,p}(\uu) w_{j}}
	  && \text{with} && \beta = \frac{ \sum_{j=0}^{n}N^\prime_{j,p}(\uu) w_{j} }{\sum_{j=0}^{n}N_{j,p}(\uu) w_{j}}
\end{align}	
For bivariate NURBS basis functions with the parametric coordinates $\uu$ and $\vv$ and the correlated orders $p$ and $q$, the first derivative in one parametric direction (i.e. $\uu$) is 
\begin{align}
	\label{eq:NURBS_Rpq_Deriv}
	\frac{\partial R_{i,j,p,q}(\uu,\vv)}{\partial \uu} 
	= \frac{
	N^\prime_{i,p}(\uu) N_{j,q}(\vv) w_{i,j}
	-  N_{i,p}(\uu) N_{j,q}(\vv) w_{i,j} \beta
	}{
	\sum_{k=0}^{n}\sum_{l=0}^{m}N_{k,p}(\uu)N_{l,q}(\vv) w_{k,l}
	}
	  && \text{with} && \beta = \frac{ \sum_{k=0}^{n}\sum_{l=0}^{m}N^\prime_{k,p}(\uu)N_{l,q}(\vv) w_{k,l} }{\sum_{k=0}^{n}\sum_{l=0}^{m}N_{k,p}(\uu)N_{l,q}(\vv) w_{k,l}}.
\end{align}


%% file: appendixOperations.tex
\section{Number of Elementary Operations}
\label{appsec:appendixOperations}

In general, elementary operations are floating point operations such as multiplication or division.
The necessary number for the evaluation of B-spline and NURBS functions is derived here.

\subsection{Basis Functions and Derivatives}

\paragraph{Univariate Basis} Using the recursive formula defined by \eqref{eq:Bspline_N0} and \eqref{eq:Bspline_Np} the computation of \revadd{all $p+1$ non-zero basis functions for a fixed parametric coordinate \protect$\uu$ results in a triangular structure of the form}{Referring to eq. B.3 directly}
\revdel{a single basis function $N_{i,p}(\uu)$ results in a triangular structure of the form\\}{Removed to condense the appendix}%
\revcomment{An equation has been removed.}%
\revdel{There are p+1 non-zero basis functions for a fixed parametric coordinate \protect$\uu$.
  Hence the previous relation extends to}{Removed to condense the appendix}%
\revcomment{An equation has been removed}%
\revdel{Fortunately, only the non-zero basis functions}{Removed to condense the appendix} %
\begin{align}
	\begin{array}{llcll}
				& 			& 		& 			& N_{i,p} 		\\
				& 			& 		& N_{i+1,p-1}	& N_{i+1,p}	\\
				& 			& \iddots	& \qquad\vdots	& \quad\vdots	\\
				& N_{i+p-1,1}	& \cdots	& N_{i+p-1,p-1}	& N_{i+p-1,p}	\\
		N_{i+p,0} 	& N_{i+p,1}	& \cdots	& N_{i+p,p-1}	& N_{i+p,p}	\\
	\end{array}.
\end{align}
\revdel{are needed. The resulting triangular}{}%
\revadd{This}{} structure consists of $p (p+1) / 2$ connections of the kind $N_{i,j-1}$ to $N_{i-1,j}$ and $N_{i,j}$.
In general, each connection needs $4$~operations (see equation~\eqref{eq:Bspline_Np}) but the result of one operation can be used twice which results in $3$~operations for each connection. 
The following total number of elementary operations is obtained:
\begin{align}
	\label{eq:op_N}
	\operations{^{N_{i}}(p)} = 3 \cdot \frac{p \left(p+1\right) }{2}
\end{align}

In order to construct all non-zero derivatives $N^\prime_{i,p}, \dots, N^\prime_{i+p,p}$ the non-zero basis functions of the previous order $p-1$ are needed first.
Then, \eqref{eq:Bspline_Np_Deriv} is applied to all non-zero terms.
The multiplication with the order $p$ can be done at once. One result of an operation can be reused again, which leads to $2$ operations for each $N^\prime_{i,p}$.
The total number of operations is
\begin{align}
	\operations{^{N^\prime_{i}}(p)} = 3 \cdot \frac{\left(p-1\right) p }{2} + 2 \left(p+1\right) 
						    = \frac{3p^2+p+4}{2}.
\end{align}
If the basis functions and their derivatives are evaluated at the same time it is
\begin{align}
	\operations{^{N_{i}+N^\prime_{i}}(p)} = 3 \cdot \frac{p \left(p+1\right) }{2} + p + 1 
						    = \frac{3p^2+5p+2}{2}.
\end{align}

For the construction of all non-zero NURBS basis functions $R_{i,p},\dots,R_{i+p,p}$ the nominator and the division is computed for each NURBS basis function. 
The dominator of~\eqref{eq:NURBS_Rpq} does not need any elementary operation, since it can be computed by a summation of each nominator.
Combined with~\eqref{eq:op_N} this results in
\begin{align}
	\operations{^{R_{i}}(p)} = 3 \cdot \frac{p \left(p+1\right) }{2} + 2 \left(p+1\right)
					   = \frac{3p^2+7p+4}{2}.
\end{align}
The computation of the first derivative of all non-zero NURBS basis functions $R^\prime_{i,p},\dots,R^\prime_{i+p,p}$ consists of the following operations: 
\begin{quote}
	\begin{tabular}{ll}
		$\left(3p^2+5p+2\right)/2$ & for all non-zero $N_{i,p}$ and $N^\prime_{i,p}$ \\
		$2 \left(p+1\right)$ & for the multiplication with their weights \\
		$1$ 			     & to set up $\beta$ of~\eqref{eq:NURBS_Rp_Deriv} \\
		$2 \left(p+1\right)$ & for the division and multiplication of~\eqref{eq:NURBS_Rp_Deriv} for each $R^\prime_{i,p}$ \\
	\end{tabular}
\end{quote}
In total, this yields to
\begin{align}
	\operations{^{R^\prime_{i}}(p)} = \frac{3p^2+13p+12}{2}.
\end{align}

\paragraph{Bivariate Basis} In order to evaluate all non-zero bivariate B-spline basis functions at a fixed parametric coordinate the following operations are need:
\begin{quote}
	\begin{tabular}{ll}
		$3 p \left(p+1\right)/2$ & for all non-zero $N_{i,p}$ \\
		$3 q \left(q+1\right)/2$  & for all non-zero $N_{j,q}$ \\
		$\left(p+1\right) \left(q+1\right)$  & for the tensor product of  $N_{i,p}$ and $N_{j,q}$ \\
	\end{tabular}
\end{quote}
Furthermore, their first partial derivatives in both parametric directions are computed by:
\begin{quote}
	\begin{tabular}{ll}
		$\left(3p^2+5p+2\right)/2$ & for all non-zero $N_{i,p}$ and $N^\prime_{i,p}$ \\
		$\left(3q^2+5q+2\right)/2$ & for all non-zero $N_{j,q}$ and $N^\prime_{j,q}$ \\
		$\left(p+1\right) \left(q+1\right)$  & for the tensor product of  $N^\prime_{i,p}$ and $N_{j,q}$ \\
		$\left(p+1\right) \left(q+1\right)$  & for the tensor product of  $N_{i,p}$ and $N^\prime_{j,q}$ \\
	\end{tabular}
\end{quote}
For the sake of clarity and simplicity, we assume that the order in each parametric direction is equal $p=q$ and obtain the following total number of elementary operations for the evaluation of bivariate B-spline basis functions $N_{i,j}$ and their derivatives in $\uu$ and $\vv$ direction $N^{{(\uu,\vv)}}_{i,j}$ at a fixed parametric coordinate:
\begin{align}
	\operations{^{N_{i,j}}(p)} &= 4 p^2 + 5 p + 1 \\
	\operations{^{N^{{(\uu,\vv)}}_{i,j}}(p)} &= 5 p^2 + 9 p + 4
\end{align}
For bivariate NURBS basis functions the number of operations is:
\begin{quote}
	\begin{tabular}{ll}
		$3 p \left(p+1\right)/2$ & for all non-zero $N_{i,p}$ \\
		$3 q \left(q+1\right)/2$  & for all non-zero $N_{j,q}$ \\
		$\left(p+1\right) \left(q+1\right)$  & for the tensor product of  $N_{i,p}$ and $N_{j,q}$ \\
		$\left(p+1\right)\left(q+1\right)$ & for the multiplication with their weights  for each $R_{i,j,p,q}$ \\
		$\left(p+1\right) \left(q+1\right) $ & division of nominator and dominator for each $R_{i,j,p,q}$\\
	\end{tabular}
\end{quote}
and for their first partial derivatives in both parametric direction:
\begin{quote}
	\begin{tabular}{ll}
		$\left(3p^2+5p+2\right)/2$ & for all non-zero $N_{i,p}$ and $N^\prime_{i,p}$ \\
		$\left(3q^2+5q+2\right)/2$ & for all non-zero $N_{j,q}$ and $N^\prime_{j,q}$ \\
		$2 \left(p+1\right) \left(q+1\right)$  & for the multiplication of $N_{j,q} w_{i,j}$ and $N^\prime_{j,q} w_{i,j}$\\
		$\left(p+1\right) \left(q+1\right)$  & for the tensor product of  $N_{i,p}$ and $N_{j,q} w_{i,j}$ \\
		$\left(p+1\right) \left(q+1\right)$  & for the tensor product of  $N^\prime_{i,p}$ and $N_{j,q} w_{i,j}$ \\		
		$\left(p+1\right) \left(q+1\right)$  & for the tensor product of  $N_{i,p}$ and $N^\prime_{j,q} w_{i,j}$ \\
		$2$ 			     & to set up $\beta$ of~\eqref{eq:NURBS_Rpq_Deriv} for $\uu$ and $\vv$ direction\\
		$4 \left(p+1\right) \left(q+1\right) $ & multiplication and devision of~\eqref{eq:NURBS_Rpq_Deriv} for $\uu$ and $\vv$ direction\\
	\end{tabular}
\end{quote}
Hence the operations for the evaluation of bivariate NURBS basis functions are
\begin{align}
	\operations{^{R_{i,j}}(p)} = 6 p^2 + 9 p + 3
\end{align}
and for their derivatives in both intrinsic directions
\begin{align}
	\operations{^{R^{(\uu,\vv)}_{i,j}}(p)} = 12 p^2 + 23 p + 13.
\end{align}

\subsection{Curve and Surface Evaluations}

The number of elementary operations for the construction of all non-zero basis functions and their derivatives are already known from the previous section.
The operations for the mapping \eqref{eq:geo_mapping} must be added to derive the number of elementary operations for the evaluation of points of B-spline and NURBS objects. 
If the order $p$ in all parametric directions is constant, the additional operations are
\begin{align}
	\label{eq:OPmapping}
	\operations{^{Map}(p)} = (p+1)^d
\end{align}
where $d$ denotes the parametric dimension.
In order to obtain the tangents of B-spline and NURBS surfaces in both intrinsic directions the mapping \eqref{eq:geo_mapping} needs to be applied twice. Hence the operations for \eqref{eq:OPmapping} are doubled. 
The combination of the different operations is summarised in~\Tabref{tab:geometryEvaluation}.

\begin{mytable}
  {\revcomment{This table has been moved to the appendix.}Number of elementary operations  $\operations{(p)}$  for the evaluation of geometry related information (i.e. points and tangents). For surfaces the order is chosen to be the same in each direction and tangents are calculated in both intrinsic directions.}%
  {geometryEvaluation}
  {ccccc}
  \mytableheader{
		    & Point   		       & Tangent  	 	      & Point   	 & Tangents         \\
	Basis    & on curve 		       & to curve 		      & on surface  & to surface}  
	B-spline & $\left(3p^2+5p+2\right)/2$  & $\left(3p^2+3p+6\right)/2$   & $5 p^2+7 p+2$    & $7 p^2+13 p+6$   \\ 
	NURBS    & $\left(3p^2+9p+6\right)/2$ & $\left(3p^2+15p+14\right)/2$ & $7 p^2+11 p+4$   & $14 p^2+27 p+15$  \\  
\end{mytable}%


%% file: integration.tex
\section{Numerical Integration}
\label{sec:integration}

The evaluation of the coefficients for the system matrices \eqref{eq:BIEdisEntries} is performed numerically. Therefore, the integration over a NURBS patch $\be$ is treated element-wise. These elements~$\ieKs$ are defined by non-zero knot spans in the knot vector of the description for the Cauchy data. Based on a subdivision scheme described next, $\ieKs$ is split into a number of integration regions $\ie$. 
As a consequence, an entry in the system matrix is a sum of integrals such as for example
\begin{equation}
  \label{eq:bem_integral}
  \mat{V}[i,j] = \sum_{e=1}^n I_{\ie_e} = \sum_{e=1}^n \, \int_{\ie_e} \fund{U}(\pt{x}_i,\pt{y}) \psi_j\ofpt{y} \dgamma{y}
\end{equation}
for the discrete single layer potential. 
The integrand includes the fundamental solution, the NURBS basis function as well as the coordinate transformation, all of which are rational functions. %
For BEM implementations, accurate integration is crucial in order to obtain correct results. Because of the possibility to use higher order basis functions as well as the fact, that the geometry is represented exactly, fewer degrees of freedom are needed as compared to conventional BEM formulations. Thus, the integration error might dominate the convergence and, as a consequence special attention is given to the numerical treatment of integrals. The following subsections describe how the integration of regular and singular integrals over NURBS patches can be performed \revadd{such that the precision of the calculated coefficients is maintained. In order to improve the efficiency of the numerical integration more sophisticated strategies should be developed.}{}


\subsection{Regular and Nearly Singular Integration}
In the presented implementation, a combination of subdivision and adaptive, successive increasing integration orders is utilised. Since NURBS patches in three dimensions can be quite extensive, in a first step all elements $\ieKs$ are subdivided so that the resulting regions $\ieAr$ have almost equal edge length in the Cartesian coordinate system. A heuristic criterion is taken to identify integrals which are \emph{nearly singular}. That is, the diameter of the integration region is compared to its distance to the collocation point $\pt{x}$. For nearly singular integration elements, a further hierarchical subdivision into regions $\ie$ is performed.
\Figref{fig:subdivision} illustrates a stretched NURBS patch $\be$ with two $\ieKs$, where in this example one is split into three $\ieAr$ and then hierarchically into several integration regions $\ie$ near the collocation point~$\pt{x}$. 
\tikzfig{tikz/subdivision}{0.9}{Subdivision of an extensive NURBS patch $\be$ where the element $\ieKs$ is subdivided into integration regions. Left: regular subdivision into $\ieAr$ with almost equal edge length and several integration regions $\ie$. Right: singular subdivision where $\ieKs_s$ denotes the area for regularisation and each $\ieKs_1$, $\ieKs_2$ and $\ieKs_3$ is treated such as a regular $\ieKs$.}{subdivision}{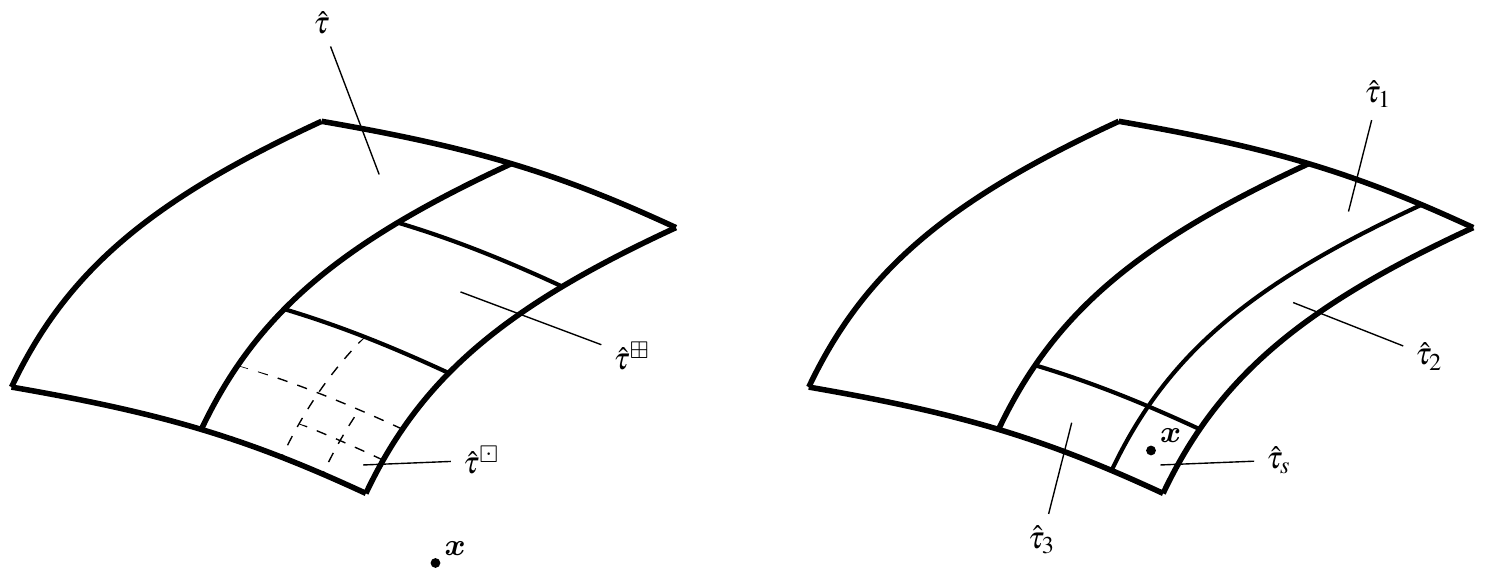}

The integral \eqref{eq:bem_integral} is carried out with standard Gauss-Legendre quadrature. Therefore, the integral on $\ie$ is transformed to an integral on a reference element $\ieg=(-1,1)^{d-1}$. The involved coordinate transformation from NURBS parameters $\pt{\uu}$ to intrinsic coordinates $\veks{\xi}=(\xi_1,\dots,\xi_{d-1})^{\trans}$ 
on $\ieg$ is defined by $\mychi_{\ieg}\ofpt{\xi}:\R^{d-1}\mapsto\R^{d-1}$. The Gauss quadrature on $\ieg$ is
\begin{equation}
  \label{eq:quadrature}
  Q = I_{\ie} + R_Q = \sum_{g=1}^{G^*} \fund{U}\left(\pt{x}_i,\pt{y}_g\right)
  \psi_j(\pt{y}_g) J_{\ieg} w_g \sqrt{g_{\ie}(\pt{\uu}_g)}
\end{equation}
with given coordinates of integration points $\pt{y}_g=\mychi(\pt{\uu}_g)=\mychi(\mychi_{\ieg}(\pt{\xi}_g))$, corresponding weights $w_g$ and Gram's determinant $g_{\ie}$ as well as the Jacobian of $\ieg$ for the integral transformation.
For each $\ie$ an adaptive integration scheme is applied. Initially, the total number of Gauss points $G^*$ is determined by the quadrature order $G_0$ in each parametric direction. Additionally, \eqref{eq:quadrature} is performed with a quadrature of higher order $G_1>G_0$. The variation of the quadrature evaluations $Q_0$ and $Q_1$ serves as a measure for the remainder 
\begin{equation}
  \label{eq:adaptive_integration}
  R_Q \sim \frac{| Q_0 - Q_1 |}{ |Q_1| } > \err_{Q,rel} 
\end{equation}
and is compared to the maximum allowed error of integration $\err_Q$. If \eqref{eq:adaptive_integration} is not satisfied, the quadrature order is increased further and the procedure restarts. 

\subsection{Singular Integration}
If the collocation point $\pt{x}$ is located on $\ieKs$, the integral is treated differently. Again, $\ieKs$ is subdivided but the procedure slightly differs from that of the regular integration. First, a rectangular area $\ieKs_s$ surrounding $\pt{x}$ is specified where the boundary integral becomes singular. This region is then subject to regularisation techniques. Duffy transformation~\cite{DUFFY1982a} is taken to treat the weakly singular integral of $\mat{V}$ in \eqref{eq:BIEdisEntries}. The procedure described by \citet{Guiggiani1990a} is used for strongly singular integrals appearing for entries in $\mat{K}$. The leftover area of integration is consequently split with the same strategy as described for regular integration. That is, equilibrating aspect ratios and hierarchical subdivision based on a geometry driven criterion. The procedure is sketched on the right in \figref{fig:subdivision}.

\revdel{The presented heuristic approach was found to be reliable when dealing with the boundary integrals which occur in the described formulation. This is confirmed by the convergence results found in \protect\secref{sec:results}.}{}%
\revdel{A further study on different and sophisticated adaptive integration strategies may greatly improve performance, but is beyond the scope of this paper.}{}%


%% file: literature.tex

\bibliographystyle{elsarticle-num-names}

\bibliography{\CommonPath/mybibliography}
